\documentclass[11pt]{article}
\usepackage{amsmath, amsthm, amsfonts, amssymb, color}
\usepackage{mathrsfs, bbm}

\allowdisplaybreaks

\newtheorem{thm}{Theorem}[section]
\newtheorem{cor}[thm]{Corollary}
\newtheorem{lem}[thm]{Lemma}
\newtheorem{remark}[thm]{Remark}
\newtheorem{prp}[thm]{Proposition}

\newtheorem{exa}[thm]{Example}
\theoremstyle{definition}
\newtheorem{defn}{Definition}[section]
\newcommand{\scr}[1]{\mathscr #1}
\definecolor{wco}{rgb}{0.5,0.2,0.3}

\numberwithin{equation}{section} \theoremstyle{remark}

\def\1{{\mathbbm 1}}

\def\wt{\widetilde}

\def\R{\mathbb R}

\def\E{\mathbb E}

\def\L{\mathcal L}

\def\F{\mathcal F}

\def\<{\langle}
\def\>{\rangle}

\def\D{\scr D}

\def\pf{\noindent{\bf Proof.} }

 \def\beq{\begin{equation}}

 \def\P{\mathbb P} 
  \def\ee{\varepsilon}

\begin{document}
\bibliographystyle{plain}

\title{\Large \bf Dirichlet heat kernel estimates of subordinate diffusion  processes with diffusive components in $C^{1, \alpha}$ open sets }

\author{
{\bf Jie-Ming Wang}
\thanks{Research partially supported  by
 NNSFC   Grant   11731009.}
}

\date{(April 29, 2024)}

\maketitle

\begin{abstract}
In this paper, we derive explicit sharp two-sided estimates of the Dirichlet heat kernels for a class of symmetric subordinate diffusion  processes with  diffusive components  in $C^{1, \alpha}(\alpha\in (0, 1])$ open sets in $\R^d$ when the scaling order  of the Laplace exponent of  purely discontinuous part of the subordinator is between $0$ and $1$ including $1.$
The main result of this paper shows the stability of Dirichlet heat kernel estimates for such processes in $C^{1, \alpha}$ open sets in the sense that the estimates  depend on the divergence  elliptic operator only via its uniform ellipticity constant  and the Dini continuity modulus of the diffusion coefficients.
As a corollary, we obtain the sharp two-sided estimates for Green functions of those processes in bounded $C^{1, \alpha}$ open sets.
\end{abstract}

\bigskip
\noindent {\bf AMS 2020 Mathematics Subject Classification}: 60J35, 47G20, 60J45, 47D07

\bigskip\noindent
{\bf Keywords and phrases}: Dirichlet heat kernel, Green function, subordinate diffusion process, second order elliptic operator, subordinate Brownian motion

\section{Introduction}
The  study of heat kernel and its estimates takes up an important place in both  analysis and probability theory.
In analysis,  heat kernel for an operator is the fundamental solution of the corresponding heat equation. When $X$ is a Markov process with the infinitesimal generator $\L_0,$ the transition density of $X$ is the fundamental solution for the operator $\L_0.$ 
For an open subset $D,$ the transition density of the subprocess of $X$ killed upon leaving $D$ (called the Dirichlet heat kernel) is the fundamental solution for the operator $\L_0$ with zero exterior condition.

Two-sided  heat kernel estimates for diffusions in $\R^d$ have a long history and many celebrated results have been established,
 see \cite{Ar, D, D2, DS}  and the references therein.
  The Dirichlet heat kernel  estimates  for the Laplace operator in $C^{1, 1}$ open sets have been established in Davies \cite{D1, D2} and Davies and Simon \cite{DS} for the upper bound estimates and in Zhang \cite{Zh} for the lower bound estimates.
Cho \cite{Cho} gives the two-sided Dirichlet heat kernel estimates for parabolic  operators of divergence form in $C^{1, \alpha}$  open sets  for
$\alpha\in (0, 1]$ under some Dini conditions on the  coefficients of the diffusion operators.

There have been intensive studies on heat kernel estimates for non-local operators in the past two decades
due to their importance in theory and applications.
 See the references \cite{BL, BBCK, BKKL, BGR3, CKK1, CKK2, CK1, CK2, CK3} therein for the  heat kernel estimates of symmetric non-local operators. For the Dirichlet heat kernel estimates for non-local operators, Chen, Kim and Song \cite{CKS1} first  established the two-sided Dirichlet heat kernel estimates for the fractional Laplacian $-(-\Delta)^{\beta/2} (\beta\in (0, 2))$ in $C^{1, 1}$ open sets in $\R^d$.
Dirichlet heat kernel estimates for more pure jump  processes  in $C^{1, 1}$ open sets in $\R^d$ have been later established, including subordinate Bownian motions, censored stable-like processes and a large class of rotationally symmetric pure jump L\'evy processes, see \cite{CKS2, CKS3, CKS5, BGR2, KM} etc..
Recently \cite{GKK} studied the Dirichlet heat kernel estimates for  symmetric  jump processes which are not necessarily L\'evy processes under some conditions on the jumping density  in $C^{1, 1}$ open sets in $\R^d.$

The boundary behavior of discontinuous processes with Gaussian components is usually  different from  pure jump processes, there are two ways for such processes to exit an open set, that is exiting continuously through the boundary or jumping across the boundary to the outside of the open set. Chen, Kim and Song \cite{CKS4, CKS7} established the Dirichlet heat kernel estimates for $\Delta+\Delta^{\beta/2}(\beta\in (0, 2))$ and a large class of subordinate Brownian motions with Gaussian components on $C^{1, 1}$ open sets when the scaling order of the pure jump part of the subordinator is strictly between $0$ and $1$. Recently,  Bae and Kim \cite{BK} extended this result to subordinate Brownian motions with Gaussian components in $C^{1, 1}$ open sets for which the scaling order of the pure jump part of the subordinator  is between $0$ and $1$ including $1$. Subordinate diffusion processes are natural extensions of subordinate Brownian motions. In this paper,  motivated by the literatures \cite{CKS4, CKS7,BK},
 we are concerned with the Dirichlet heat kernel estimates for symmetric subordinate diffusions with diffusive components when the scaling order of pure discontinuous part of the subordinator is between $0$ and $1$ including $1$ in
$C^{1, \alpha}$ open sets in $\R^d$  for $\alpha\in (0, 1].$

More specifically,  we consider the following second order elliptic operator   in divergence form in $\R^d (d\geq 1)$
\begin{equation}\label{e:1.1}
\L_0 f(x) =\dfrac{1}{2}\sum_{i, j=1}^d \dfrac{\partial}{\partial x_i}\left(a_{ij}(x)\dfrac{\partial}{\partial x_j} f(x)\right),
\end{equation}
where $(a_{ij}(x))_{1\leq  i, j\leq d}$ is a symmetric $d\times d$ matrix-valued function on $\R^d$
that is   uniformly bounded and elliptic;  that is,
there exists a constant
 $\lambda_0\geq 1$
such that for all $x\in\R^d$ and $\xi\in\R^d,$
\begin{equation}\label{e:1.2}
\lambda_0^{-1} |\xi|^2 \leq  \sum_{i, j=1}^d a_{ij}(x) \xi_i \xi_j \leq \lambda_0  | \xi|^2.
\end{equation}
Associated with $\L_0$ is a symmetric diffusion process $X$ in $\R^d$ whose associated Dirichlet form $(\mathcal{E}_X, \F)$ on $(L^2(\R^d); dx)$ is given by
$$\mathcal{E}_X(u, v)=\dfrac{1}{2}\int_{\R^d} A(x)\nabla u(x)\cdot \nabla v(x) \,dx, \quad \F=\overline{C_c^1(\R^d)}^{\mathcal{E}_1},$$
where $A(x)=(a_{ij}(x))_{1\leq  i, j\leq d}$ and $\mathcal{E}_1(u, v):=\mathcal{E}_X(u, v)+\int_{\R^d} u(x)v(x)\,dx.$
It is well-known that $X$ has a jointly H\"older continuous transition density function $p^X(t, x, y)$, which enjoys the following celebrated Aronson’s two-sided heat kernel estimates: there are constants $c_k=c_k(d, \lambda_0) > 0, k =1,...,4,$ so that
\begin{equation}\label{e:1.3'}
c_1 p^W(t, c_2x, c_2y)\leq p^X(t, x, y)\leq c_3p^W(t, c_4x, c_4y) \quad \mbox{for} \quad t>0 \quad \mbox{and} \: x,y\in\R^d,
\end{equation}
where $p^W(t, x, y)$ is the transition density function of a Brownian motion in $\R^d.$

A subordinator $S_t$ is an increasing L\'evy process in $\R_{+}$ starting from $0$,  which can be characterized through its Laplace exponent $\varphi:$
$\E[e^{-\lambda S_t}]=e^{-t\varphi(\lambda)}$, $\lambda\geq 0.$ The Laplace exponent of a subordinator belongs to the class of  Berstein functions
$\mathcal{BF}=\{f\in C^\infty(0, \infty): f\geq 0,  (-1)^{n-1}f^{(n)}\geq 0, n\in \mathbb{N}\}$ and  has the representation
$$
 \varphi (\lambda)=b\lambda+\phi(\lambda) \quad {\rm with} \quad \phi(\lambda):=\int_0^\infty (1-e^{-\lambda t})\,\mu(dt),
$$
where $\phi$ is the Laplace exponent of the pure jump part of the subordinator $S_t,$  $b\geq 0$ is called the drift of $S$ and $\mu$ is a  measure (called the L\'evy measure of  $\varphi$)
 on $[0, \infty)$ satisfying $\int_0^\infty (1\wedge t)\mu(dt)<\infty$.
 The Laplace exponent $\varphi$ of a subordinator $S$ is said to be  a complete Bernstein function
if its L\'evy measure $\mu (dt)$
 has a completely monotone density $\mu(t)$ with respect to the Lebesgue measure on $(0, \infty)$;
that is, if $\mu (dt) = \mu  (t) dt$ with $\mu \in C^\infty (0, \infty)$ and $(-1)^n \mu ^{(n)}(t) \geq 0$ on $(0, \infty)$ for every integer $n\geq 0$.
In this case, we say $S$ is a complete subordinator.
It is known that most familiar Bernstein functions are complete Bernstein functions.
See \cite[Chapter 15]{SSV} for  more details of complete Bernstein functions.

Throughout this paper, we let $Y_t:=X_{S_t}$ be a diffusion $X$ subordinated by a subordinator $S$, where  $S$ is a complete subordinator independent of $X$ with the positive drift $b>0.$ Without loss of generality, we assume $b=1.$
The process $Y$ is a symmetric Hunt process with the transition density function
\begin{equation}\label{e:1.5n}
p(t, x, y)=\int_0^\infty p^X(s, x, y) \P(S_t\in ds).
\end{equation}
The Dirichlet form $(\mathcal{E}, \D(\mathcal{E}))$ of $Y$ on $(L^2(\R^d); dx)$ is given by
\begin{equation}\label{e:1.6'}\begin{aligned}
\mathcal{E}(u, v)=\dfrac{1}{2}&\int_{\R^d} A(x)\nabla u(x)\cdot \nabla v(x) \,dx\\
&+\int_{\R^d\times\R^d} (u(x)-u(y))(v(x)-v(y)) J(x, y)\,dx\,dy
\end{aligned}\end{equation}
and $\D(\mathcal{E})=\D(\mathcal{E}_X)$ (see \cite[Theorem 2.1]{Ok}),
where
\begin{equation}\label{e:1.6n}
J(x, y)=\int_0^\infty p^X(t, x, y)\mu(t)\,dt.
\end{equation}
It is known that any Hunt process admits a L\'evy system that describes how the process jumps.
By applying the  similar argument in \cite[Lemma 4.7]{CK1}, for any  nonnegative function $f$ on $\R_+\times
\R^d\times \R^d$ vanishing along the diagonal of $\R^d \times \R^d$, for any stopping time $T$
with respect to the minimal admissible  augmented  filtration generated by $Y$ and $x\in \R^d$,
\begin{equation}\label{e:1.6}
\E_x \left[ \sum_{s\leq T}f(s, Y_{s-}, Y_s); Y_{s-}\neq Y_s \right]
=\E_x \left[ \int_0^T\int_{\R^d}
f(s,Y_s,y) J(Y_{s-}, y)\,dy\,ds \right].
\end{equation}

We introduce the following scaling conditions for a function $g: (0, \infty)$ to $(0, \infty).$

\begin{defn}
Suppose  $g$ is a function from $(0, \infty)$ to $(0, \infty).$

  (1)  We say that $g$ satisfies  $L^a(\gamma, c_L)$ (resp. $L_a(\gamma, c_L)$) if there exist $a\geq 0, \gamma>0$ and $c_L\in (0, 1]$ such that
$$\dfrac{g(R)}{g(r)}\geq c_L (\frac{R}{r})^\gamma \quad \mbox{for all}\quad  a< r\leq R \quad (\mbox{resp.}\: 0<r\leq R\leq a).$$

(2) We say that $g$ satisfies  $U^a(\delta, C_U)$ (resp. $U_a(\delta, C_U)$) if there exist $a\geq 0, \delta>0$ and $C_U\in [1, \infty)$ such that
$$\dfrac{g(R)}{g(r)}\leq C_U (\frac{R}{r})^\delta \quad \mbox{for all} \quad a< r\leq R \quad (\mbox{resp.}\: 0<r\leq R\leq a).$$

\end{defn}

We define
 $$H(\lambda):=\phi(\lambda)-\lambda\phi'(\lambda).$$
 The function $H$ appeared early in the work of Jain and Pruitt \cite{JP}.
 When the scaling order of the Laplace exponent of the pure jump part of the subordinator  is not strictly less than $1$, the heat kernel estimates of the subordinate Brownian motion will have a new form which is associated closely with the function $H$;  see \cite{M} and \cite{BK}.

Let $W$ be a Brownian motion in $\R^d$ independent of $S_t$ and denote by $Y^0_t:=W_{S_t}$ the Brownian motion $W$ subordinated by the subordinator $S_t.$
The two-sided heat kernel estimates for $Y^0$ in $\R^d$ have been established in \cite[Theorem 1.3]{BK} under some mild conditions on $H$ in Theorem \ref{T0} below.
Denote by $p^0(t, x, y)$  the transition density function of $Y^0$ in $\R^d.$
In view of \eqref{e:1.3'} and \eqref{e:1.5n}, it is easy to see that there are constants $c_k=c_k(d, \lambda_0) > 0, k =1,...,4$ so that
\begin{equation}\label{e:1.7}
c_1 p^0(t, c_2x, c_2y)\leq p(t, x, y)\leq c_3p^0(t, c_4x, c_4y) \quad \mbox{for} \quad t>0 \quad \mbox{and} \: x,y\in\R^d,
\end{equation}
By \eqref{e:1.7} and the result in  \cite[Theorem 1.3]{BK}, we have the following result.

\begin{thm}\label{T0}
Suppose $a_{ij}(x)$ satisfies the uniformly elliptic condition \eqref{e:1.2}.
Suppose $H$ satisfies $L^a(\gamma, c_L)$ and $U^a(\delta, C_U)$ for some $a>0$ with  $\delta<2,$ then for each $T>0$ and $M>0,$ there exist
 $C_1=C_1(d, \lambda_0,  \phi, M, T)$ and $C_k=C_k(d, \lambda_0,  \phi), k=2, \cdots, 5$ such that for any $t\in (0, T]$ and $|x-y|\leq M/2,$
\begin{equation}\label{e:2.26'}
C_1t^{-d/2}\wedge\left(p^W(t, C_2x, C_2y)+q(t, C_3x, C_3y)\right)\leq p(t, x, y)\leq  C_1t^{-d/2}\wedge\left(p^W(t, C_4x, C_4y)+q(t, C_5x, C_5y)\right),
\end{equation}
where
$$q(t, x, y)=\dfrac{tH(|x-y|^{-2})}{|x-y|^d}+\phi^{-1}(1/t)^{d/2}e^{-|x-y|^2\phi^{-1}(1/t)}.$$
If $H$ satisfies $L^0(\gamma, c_L)$ and $U^0(\delta, C_U)$ with  $\delta< 2,$ then there exist $C_k=C_k(d, \lambda_0,  \phi), k=6, \cdots, 10$ such that for any $t>0$ and $x, y\in\R^d,$
\begin{equation}\label{e:2.29'}
\begin{aligned}
&C_6^{-1}(t^{-d/2}\wedge (\phi^{-1}(1/t))^{d/2})\wedge\left(p^W(t, C_7x, C_7y)+q(t, C_8x, C_8y)\right)\leq p(t, x, y)\\
&\quad\quad\leq  C_6(t^{-d/2}\wedge (\phi^{-1}(1/t))^{d/2})\wedge\left(p^W(t, C_9x, C_9y)+q(t, C_{10}x, C_{10}y)\right).
\end{aligned}
\end{equation}
\end{thm}

\begin{remark} \rm
$H(\lambda)$ and $\phi(\lambda)$ are comparable for $\lambda>a$ if and only if $\phi$ satisﬁes $U^a(\delta, C_U)$ with $\delta<1$ for some $a\geq 0$(see \cite[Proposition 2.9]{M}).
That is,  $H$  can be replaced by  $\phi$ in Theorem \ref{T0} when the upper scaling order $\delta$ of $\phi$ is strictly less than $1.$
Observing that $J(x, y)$ satisfies the
conditions in \cite{CK3} under the condition that  $\phi$ satisfies $L^0(\gamma, c_L)$ and $U^0(\delta, C_U)$ with  $\delta< 1,$ thus $q(t, x, y)$ in Theorem \ref{T0} is indeed comparable to $t\phi(|x-y|^{-2})/|x-y|^d$ by \cite{CK3} in this case.

 Although the dependence of the multiplying constants in \cite[Theorem 1.3]{BK} are implicit, while by carefully checking the proof and by \eqref{e:1.7}, the dependence of the constants in Theorem \ref{T0} can be obtained.
\end{remark}

\smallskip
Recall that an open set $D$ in $\R^d$ (when $d\geq 2$) is said to be $C^{1,\alpha}(\alpha\in (0, 1])$ if there exist a localization radius
$R_0 >0$ and a constant $\Lambda_0>0$ such that for every $z\in\partial D,$ there exist a $C^{1,\alpha}$ function
$\Gamma=\Gamma_z: \R^{d-1}\rightarrow \R$ satisfying
\begin{equation}\label{e:1.5}
\Gamma(0)=\nabla \Gamma(0)=0, \|\nabla \Gamma\|_\infty\leq \Lambda_0 ,
|\nabla \Gamma(x)-\nabla \Gamma(y)|\leq \Lambda_0 |x-y|^\alpha,
\end{equation}
 and an orthonormal coordinate system
$CS_z: y=(y_1, \cdots, y_{d-1}, y_d)=:(\wt {y}, y_d)\in \R^{d-1}\times \R$ with its origin at $z$ such that
$$
B(z, R_0)\cap D=\{y=(\wt {y}, y_d)\in B(0,R_0) \hbox{ in } CS_z: y_d>\Gamma(\wt {y})\}.
$$
The pair $(R_0, \Lambda_0)$ is called the characteristics of the $C^{1,\alpha}$ open set $D.$
Without loss of generality, throughout this paper, we assume that the characteristics $(R_0, \Lambda_0)$
of a $C^{1,\alpha}$ open set satisfies $R_0\leq 1$ and $\Lambda_0\geq 1$.
For any $x\in D,$ let $\delta_D(x)$ denote the Euclidean distance between $x$ and $D^c$.

We say that the path distance in a domain (connected open set) $U$ is comparable to the Euclidean distance with characteristic $\chi_1$
if for every $x$ and $y$ in $U,$ there is a rectifiable curve $l$ in $U$ which connects $x$ to $y$ such that the length of $l$ is less than or equal to $\chi_1|x-y|.$ Clearly such a property holds for all bounded $C^{1, \alpha}$ domains, $C^{1, \alpha}$ domains with compact complements and domain consisting of all the points above the graph of  $C^{1, \alpha}$ function.

\medskip

To establish the Dirichlet heat kernel estimates for the process $Y$ in  $C^{1, \alpha}$ open sets,
we need some additional conditions. We assume  the entries $a_{ij}(x)$, $1\leq i, j\leq d$,  are Dini continuous, that is,
\begin{equation}\label{e:1.3}
\sum_{i, j=1}^d|a_{ij}(x )-a_{ij}(y)|\leq \ell (|x - y|) \quad \hbox{for all } x ,  y\in\R^d \hbox{ and } 1\leq i, j\leq d,
\end{equation}
 where  $\ell (\cdot): [0, \infty)\rightarrow [0, \infty)$ is   an increasing
 continuous function  with $\ell (0)=0$ and $\int_0^1 \ell (t)/t\, dt<\infty.$
Let $j(x, y)$ be the jumping kernel of the subordinate Brownian motion $Y^0_t=W_{S_t}.$
Since $W_{S_t}$ is rotationally symmetric, we also write $j(x, y)=j(|x-y|).$
We have
$$j(x, y)=j(|x-y|)=\int_0^\infty p^W(t, x, y)\mu(t)dt, \quad x, y\in\R^d,$$
It follows from\eqref{e:1.3'} and \eqref{e:1.6n} that there exist $c_k=c_k(d, \lambda_0)>1, k=1, 2, 3$ such that
\begin{equation}\label{e:2.8'}
c_1^{-1} j(c_2|x-y|)\leq J(x, y)\leq c_1j(c_3|x-y|), \quad x, y\in \R^d.
\end{equation}
By Theorem 2.1 and Lemma 2.4 in \cite{BK}, when $H$ satisfies $L^a(\gamma, c_L)$ for some $a>0,$
there exists a constant $c$  such that
 \begin{equation}\label{e:2.17}
  j(|x-y|)\leq c\dfrac{H(|x-y|^{-2})}{|x-y|^d}, \quad x, y\in\R^d\setminus {\rm diag}.
 \end{equation}
In particular, when $H$ satisfies $L^a(\gamma, c_L)$ and $U^a(\delta, C_U)$ with  $\delta<2$ for some $a>0,$
 for each $M>0,$ there exists a positive constant $c$ depending on $M$ such that
 \begin{equation}\label{e:2.16}
 c^{-1}\dfrac{H(|x-y|^{-2})}{|x-y|^d} \leq j(|x-y|)\leq c\dfrac{H(|x-y|^{-2})}{|x-y|^d}, \quad x, y\in\R^d\setminus {\rm diag}\quad {\rm with} \quad |x-y|\leq M.
 \end{equation}
Furthermore, when $H$ satisfies $L^0(\gamma, c_L)$ and $U^0(\delta, C_U)$ with  $\delta<2,$
  there exists a positive constant $c$ such that
 \begin{equation}\label{e:1.15}
 c^{-1}\dfrac{H(|x-y|^{-2})}{|x-y|^d} \leq j(|x-y|)\leq c\dfrac{H(|x-y|^{-2})}{|x-y|^d}, \quad x, y\in\R^d\setminus {\rm diag}.
 \end{equation}
 Since $Y^0_t=W_{S_t}$ is a L\'evy process in $\R^d$, the jumping density function $j(r)$  satisfies  $\int_{\R^d} (1\wedge |z|^2) j(|z|)\,dz<\infty.$
Note that the jumping density function $j(r)$ is non-increasing, it follows from \cite[(1.5)]{CKS7} that there exists $c>0$ such that $j(r)\leq cr^{-(d+2)}$ for $r\in (0, 1).$
 When $H$ satisfies $L^a(\gamma, c_L)$ and $U^a(\delta, C_U)$ with  $\delta<2$ for some $a>0,$ in view of \eqref{e:2.16}, $j(r)\asymp \frac{H(r^{-2})}{r^d}$ for $r\in (0, 1),$ thus there exists $c>0$ such that
$H(\lambda)\leq c\lambda$ for $\lambda>1.$

\medskip
 In the following, we consider the following assumptions on $H.$

\begin{enumerate}
\item[{\bf (A1)}]  $H$ satisfies $L^a(\gamma, c_L)$ and $U^a(\delta, C_U)$ with  $\delta< 1$ for some $a>0,$
or  $H$ satisfies $L^a(\gamma, c_L)$ and $U^a(\delta, C_U)$ with  $\delta=  1$ for some $a>0$ and $\gamma>1/2.$

\item[{\bf (A2)}] $H$ satisfies $L^0(\gamma, c_L)$ and  $U^0(\delta, C_U)$ with $\delta< 2.$

 \end{enumerate}

 The assumption (A1)  shows the  conditions on the scaling order of $H$ near the infinity.
Note that $H(\lambda)\leq c\lambda$ for $\lambda>1.$
 The condition $\delta\leq 1$ in the assumption (A1) is in fact a mild condition.
Note that $H(\lambda)$ and $\phi(\lambda)$ are comparable for $\lambda>a$ if and only if $\phi$ satisﬁes $U^a(\delta, C_U)$ with $\delta<1$ for some $a\geq 0$(see \cite[Proposition 2.9]{M}),
the first condition in the assumption (A1)  is equivalent that  the scaling order of $\phi$ near the infinity is strictly between $0$ and $1.$
For the latter  of the assumption (A1) that $H$ satisfies $L^a(\gamma, c_L)$ and $U^a(\delta, C_U)$ with  $\delta=  1$ for some $a>0$ and $\gamma>1/2,$ this condition covers the case $H(\lambda)=\lambda\ell_0(\lambda)$ for $\lambda\geq 1$, which corresponds to $j(r)\asymp\frac{\ell_0(r^{-2})}{r^{d+2}}$ for $r\leq 1$ by \eqref{e:2.16}, where $\ell_0$ slowly  varies  at inﬁnity,
i.e.  $\lim_{\lambda\rightarrow\infty} \frac{\ell_0(\lambda s)}{\ell_0(\lambda)}=1$ for each $s>0.$

The assumption (A2) shows  the upper scaling order  $\delta$ of $H$ near  $0$ and the infinity is less than $2.$
Note that $j(\lambda)\asymp \frac{H(\lambda^{-2})}{\lambda^d}$ for $\lambda>0$ by \eqref{e:1.15} under the assumption (A2).
This assumption  covers the case $j(\lambda)\asymp \lambda^{-(d+\ee)}$ for $\lambda\geq 1$ with $\ee \in (0, 2\delta]$ near the infinity.

For any open set $D\subset\R^d$ and positive constants $c_1$ and $c_2,$ we define
$$\begin{aligned}
&h_{D, c_1, c_2}(t, x, y):=\left(1\wedge \dfrac{\delta_{D}(x)}{\sqrt t}\right)\left(1\wedge \dfrac{\delta_{D}(y)}{\sqrt t}\right) \\
 &\quad \times \left[t^{-d/2}\wedge\left(t^{-d/2}e^{-c_1|x-y|^2/t}+\dfrac{tH(|x-y|^{-2})}{|x-y|^d}
+\phi^{-1}(1/t)^{d/2}e^{-c_2|x-y|^2\phi^{-1}(1/t)}\right)\right].
\end{aligned}$$

The following is the main theorem of this paper. 
Let $Y^D$ be the subprocess of $Y$ killed upon leaving $D.$
 Denote by $p_D(t, x, y)$ the transition density function of $Y^D.$

\begin{thm}\label{T1}
Suppose
 $(a_{ij}(x))_{1\leq i, j\leq d}$ satisfies the conditions \eqref{e:1.2} and \eqref{e:1.3}.
Suppose that $D$ is a $C^{1, \alpha} (\alpha\in (0, 1])$ open set in $\R^d (d\geq 1)$ with characteristics $(R_0, \Lambda_0).$
If $D$ is bounded,  assume that $H$ satisfies the assumption (A1).
If $D$ is unbounded,  assume that $H$ satisfies the assumptions (A1) and (A2).

(i) For every $T>0,$ there exist positive constants $C_1=C_1(d, \lambda_0, \ell, \phi, R_0, \Lambda_0, T)$ and $a_U, b_U$ depending  on
$(d, \lambda_0,  \phi)$ such that
for  any $x, y\in D$ and $t\in (0, T),$
\begin{equation}\label{e:1.8}
p_D(t,x, y)\leq C_1h_{D, b_U, a_U}(t, x, y).
\end{equation}

(ii) Assume the path distance in each connected component of $D$ is comparable to the Euclidean distance with characteristic $\chi_1.$ For every $T>0,$ there exist positive constants $C_2=C_2(d, \lambda_0, \ell, \phi, R_0, \Lambda_0, \chi_1, T)>0$ and $a_L,  b_L$ depending  on
$(d, \lambda_0, \ell, \phi, R_0, \Lambda_0, \chi_1)$ such that
for  any $x, y\in D$ and $t\in (0, T),$
\begin{equation}\label{e:1.9}
p_D(t,x, y)\geq C_2h_{D, b_L, a_L}(t, x, y).
\end{equation}

(iii) If $D$ is bounded, then for each $T>0,$ there exists $C_3=C_3(d, \lambda_0, \ell, \phi, R_0, \Lambda_0, {\rm diam}(D), T)>1$
such that for any $(t, x, y)\in (T, \infty)\times D\times D,$
\begin{equation}\label{e:1.10}
C_3^{-1}e^{-\lambda_1 t}\delta_D(x)\delta_D(y)\leq p_D(t, x, y)\leq C_3e^{-\lambda_1 t}\delta_D(x)\delta_D(y),
\end{equation}
where $-\lambda_1<0$ is the largest eigenvalue of the generator of $Y^D.$
\end{thm}

\medskip
 The  Dirichlet heat kernel estimates for  subordinate Brownian motions with Gaussian components  in a $C^{1, 1}$ open set $D$   are obtained  in \cite{BK}  under the conditions that $H$ satisfies $L^a(\gamma, c_L)$ and $U^a(\delta, C_U)$ with  $\delta< 2$ for some $a>0$ (resp. $a=0$) when $D$ is bounded (resp. unbounded). Although there are some differences between the assumption (A1)  and the condition in \cite{BK} when $D$ is bounded, the assumption (A1) is mild and covers a large class of  subordinators. In fact, as we mentioned before, assumption (A1) covers the case when the  scaling order  of $\phi$ near the infinity  is strictly between $0$ and $1.$ If the upper scaling order $\delta$ of $\phi$ near the infinity is equal to $1,$ the second part of the assumption (A1) contains the case $H(\lambda)=\lambda\ell_0(\lambda)$ for $\lambda\geq 1,$ which is equivalent that $j(r)\asymp\frac{\ell_0(r^{-2})}{r^{d+2}}$ for $r\leq 1$ by \eqref{e:2.16}, where $\ell_0$  slowly varies  at inﬁnity. This condition  includes a large class of subordinators with Laplace exponents that vary regularly at infinity with index $1.$
 When $D$ is unbounded,  the assumptions (A1) and (A2) holding simultaneously is equivalent that $H$ satisfies (A1) near the infinity and the scaling order of $H$ near $0$ is less than $2.$
 The following are some  examples  of Theorem \ref{T1}.

\begin{exa}\label{E1}
(1) Let $\phi(\lambda)=\frac{\lambda}{\log(1+\lambda^{\beta/2})} (\beta \in (0, 2))$ be  the Laplace exponent of the conjugate gemoetric stable subordinator without killing. Then $\phi$ is a complete Bernstein function by \cite[Example 5.11 and (ii) Page 90]{BBKR}. We have
$$\phi^{-1}(\lambda) \asymp \left\{\begin{array}{ll}
\lambda^{2/(2-\beta)}, & \quad 0<\lambda<2\\
\lambda\log\lambda,  & \quad \lambda\geq 2\\
 \end{array}\right.
 \quad
 H(\lambda) \asymp \left\{\begin{array}{ll}
\lambda^{1-\beta/2}, & \quad 0<\lambda<2\\
\lambda/(\log\lambda)^2,  & \quad \lambda\geq 2.\\
 \end{array}\right.
 $$
 Then $H$ satisfies $L^0(\gamma, c_L)$ and $U^0(\delta, C_U)$ with $\delta\leq 1$ and satisfies $L^2(\gamma, c_L)$ and $U^2(\delta, C_U)$ with $\delta= 1$ and $\gamma>1/2.$ By Theorem \ref{T1}, the upper bounds and the lower bounds in \eqref{e:1.8} and \eqref{e:1.9} hold for the process $Y$ when $D$ is a (possibly unbounded) $C^{1, \alpha}$ open set.

 (2) Let $\phi(\lambda)=\frac{\lambda}{\log(1+\lambda)}-1$ be  the Laplace exponent of the conjugate gamma subordinator without killing.
Then $\phi$ is a complete Bernstein function by \cite[Example 5.10 and (ii) Page 90]{BBKR}.   We have
$$\phi^{-1}(\lambda) \asymp \left\{\begin{array}{ll}
\lambda, & \quad 0<\lambda<2\\
\lambda\log\lambda,  & \quad \lambda\geq 2\\
 \end{array}\right.
 \quad
 H(\lambda) \asymp \left\{\begin{array}{ll}
\lambda^2, & \quad 0<\lambda<2\\
\lambda/(\log\lambda)^2,  & \quad \lambda\geq 2.\\
 \end{array}\right.
 $$
 Then $H$ satisfies $L^2(\gamma, c_L)$ and $U^2(\delta, C_U)$ with $\delta=1$ and $\gamma>1/2.$
 By Theorem \ref{T1},  \eqref{e:1.8}-\eqref{e:1.10} holds for the process $Y$ when $D$ is a bounded $C^{1, \alpha}$ open set.
\end{exa}

 \medskip

 Theorem \ref{T1} shows the stability of Dirichlet heat kernel estimates for subordinate diffusion processes with diffusive  part in a $C^{1, \alpha}$ open set in the sense that the constants in Theorem \ref{T1} depend on $\L_0$ only via the uniform elliptic constant $\lambda_0$ and the Dini continuity modulus $\ell$ of the diffusion coefficients $a_{ij}(x).$
 Theorem \ref{T1} is  new for subordinate Brownian motions with Gaussian components in a less smooth $C^{1, \alpha}$ open set with $\alpha\in (0, 1)$.
The Dini condition on $a_{ij}$ in Theorem \ref{T1} is in fact a mild condition.
This condition is used in \cite{GW} for the upper bound estimates of the Green function for the divergence form second order elliptic operator  in the domain satisfying the exterior sphere condition  and in \cite{Cho} for the two-sided Dirichlet heat kernel estimates for the parabolic operator in divergence form in $C^{1, \alpha}$ open sets.

  The key ingredient   in the previous literatures \cite{CKS4, CKS7,BK} for the Dirichlet heat kernel estimates for subordinate Brownian motions with Gaussian components in $C^{1, 1}$ open sets is   the test function method and Dynkin's formula.
  Since the space of smooth functions with compact support are contained in the domains of the inﬁnitesimal generator of the L\'evy process, by choosing appropriate test functions and  computing the generator acting on  the test function  and using Dynkin's formula, the exit time and the exit distribution estimates  for the subordinate Brownian motions with Gaussian components  could be  obtained.
 Based on these estimates, the decay rate of the Dirichlet heat kernel near the boundary can be established.
 While in our case, as the smooth function with compact support may not be contained in the domains of the infinitesimal  generators of the process $Y$ in this paper and  the process $Y$  may not be a semimartingale, this makes it difficult to adapt the methods in the previous literatures for L\'evy processes to our case and thus causes    difficulty in our setting.

 In this paper, we mainly use probabilistic method. Instead of test function method, we  make use of the  resurrection formula between the  killed subordinate Markov process and the subordinate killed Markov process on an open set established by Song and  Vondra\v{c}ek \cite{SV1} to compute the exit time estimates  of $Y$ from a small $C^{1, \alpha}$ domain. Then we use the "box" method developed in Bass and Burdzy \cite{Ba, BB1}  to obtain the exit distribution estimates for $Y$ from a $C^{1, \alpha}$ open set. By  combining these results and following the probabilistic strategies   in \cite{CKS7, BK},  the  Dirichlet heat kernel estimates for $Y$ in a $C^{1, \alpha}$ open set can be obtained. In fact, by virtue of the resurrection formula in \cite{SV1}, we prove the exit time of $Y$ from a  $C^{1, \alpha}$ domain can be written as the sum in terms of the Green function of subordinate killed diffusion $Z^D=X^D(S_t)$ and the resurrection kernel $q_D(y, z)$ (see \eqref{e:2.21'} and \eqref{e:2.22} below) and then use this formulation to obtain the estimates of the exit time of $Y$. This method may be used for the study of more classes of subordinate Markov processes.

\medskip
Define
$$
g_D(x, y):=\left\{\begin{array}{ll}
\left(1\wedge \dfrac{\delta_D(x) \delta_D(y)}{|x-y|^2}\right)|x-y|^{2-d}, & \quad d\geq 3\\
\log\left(1+\dfrac{\delta_D(x)\delta_D(y)}{|x-y|^2}\right),  & \quad d=2\\
 (\delta_D(x)\delta_D(y))^{1/2} \wedge \dfrac{\delta_D(x)\delta_D(y)}{|x-y|}, & \quad d=1.
\end{array}\right.
$$

It is known from \cite[Theorem 4.8]{CKP} that $g_D(x, y)$ is comparable to the Green function of the diffusion $X$ on bounded $C^{1, \alpha}$ domains (i.e. connected open sets) in $\R^d$ under the conditions \eqref{e:1.2} and \eqref{e:1.3}. By integrating the two-sided heat kernel estimates in Theorem \ref{T1} with respect to $t,$ we can obtain two-sided estimates of Green function
$G_D(x, y):=\int_0^\infty p_D(t, x, y)\,dt$ of $Y$ for a bounded  $C^{1, \alpha}$ open set  in $\R^d$.

\begin{cor}\label{C1}
Suppose $(a_{ij}(x))_{1\leq i, j\leq d}$ satisfies the conditions \eqref{e:1.2} and \eqref{e:1.3}.
Suppose that $D$ is a bounded $C^{1, \alpha} (\alpha\in (0, 1])$ open set in $\R^d$ with characteristics $(R_0, \Lambda_0)$
and $H$ satisfies   the assumption (A1). Then
there exists $C=C(d, \lambda_0, \ell, \phi, R_0, \Lambda_0, {\rm diam}(D))>1$ such that
$$C^{-1}g_D(x, y)\leq G_D(x, y)\leq Cg_D(x, y), \quad x, y\in D.$$
\end{cor}

The organization of this paper is as follows. In Section 2, by using the resurrection formula between the  subordinate killed Markov process and the killed subordinate  Markov process on an open set in \cite{SV1}, we  obtain  the two-sided estimates of the  exit time for $Y$ from a small bounded $C^{1, \alpha}$ domain is  comparable to the exit time of a Brownian motion. In Section 3, we use the exit time estimates in Section 2 and  the "box" method developed in Bass  and Burdzy \cite{Ba, BB1} to establish  the exit distribution estimates for the process $Y$ from  a $C^{1, \alpha}$ open set $D$. In Section 4, by applying the routine argument in \cite{CKS7, BK} and the results in Sections 2 and 3, we establish the two-sided  Dirichlet heat kernel estimates of $Y$ in a $C^{1, \alpha}$ open set $D$.

\section{ Exit time estimates from a small bounded $C^{1, \alpha}$ domain}

Throughout this paper, unless specified we assume $d\geq 1$. Let $X$ be a diffusion process associated with  $\L_0$ under the conditions \eqref{e:1.2} and \eqref{e:1.3}.
Let  $S$ be a complete subordinator  independent of $X$ with the Laplace exponent given by
 $$\E \exp(-\lambda S_t)=\exp(-t\varphi(\lambda)), \quad \lambda>0,$$
where the Laplace exponent
$$\varphi(\lambda)=\lambda+\int_0^\infty (1-e^{-\lambda t})\mu(t)\,dt,$$
where $\mu \in C^\infty (0, \infty)$ and $(-1)^n \mu ^{(n)}(t) \geq 0$ on $(0, \infty)$ for every integer $n\geq 0$.
Let $$Y_t:=X_{S_t}.$$
In this section, we shall use the resurrection formula from \cite{SV1} between the killed subordinate Markov process and the subordinate killed Makov process to derive the estimates of the exiting time of $Y$ from a small domain $rD$ with $r\in (0, 1),$
where $D$ is  bounded $C^{1, \alpha}$ domain with characteristics $(R_0, \Lambda_0).$

 For an open set $D,$ let $X^D$ be the part process of $X$ killed upon leaving $D$ and let $Z^{D}_t:=X^{D}(S_t)$.
The process $Z^D$ is called  the subordinate killed diffusion in $D$.
 We will use $\zeta$ to denote the life time of the process $Z^D_t.$
 It follows from  \cite[(4.2)]{SV1} that the subordinate killed process $Z^D$ admits a L\'evy system of the form
 $(J^{Z^D}(x, y)\,dy, dt),$ where
 \begin{equation}\label{e:2.2}
 J^{Z^D}(x, y):=\int_0^\infty p^X_D(t, x, y) \,\mu(t)dt, \quad x, y\in D.
 \end{equation}
 Let $\{\F_t\}_{t\geq 0}$ be the usual argumentation of the natural filtration generated by the diffusion process $X.$
 Let $\tau$ be an $(\F_{t+})$-stopping time. Define
$$
 \sigma_\tau:=\inf\{t>0: S_t>\tau\}.
$$
Let $\tau_D:=\inf\{t>0: Y_t\notin D\}$ and $\tau^X_D:=\inf\{t>0: X_t\notin D\}.$ The process $Z^D_t$ can be written as
 $$
 Z^D_t=\left\{\begin{array}{cc}
                Y_t, & t<\sigma_{\tau^X_D} \\
                \partial, & t\geq \sigma_{\tau^X_D}
              \end{array}
 \right.
 =\left\{\begin{array}{cc}
                X_{S_t}, & S_t<\tau^X_D \\
                \partial, & S_t\geq \tau^X_D.
              \end{array}
 \right.
 $$
 Let $Y^D$ be the process $Y$ killed upon leaving $D.$
This shows  that $\sigma_{\tau^X_D}\leq \tau_D$ and the process $Z^D$ is a subprocess of $Y^D$ by killing $Y^D$ at the terminal time $\sigma_{\tau^X_D}.$
By Proposition 3.2 in \cite{SV1}, the process $Y^D$ can be obtained from $Z^D$ by resurrecting the latter at most countably many times.

The potential measure of the subordinator $S_t$ is defined to be
$$U(A):=\E \int_0^\infty 1_{\{S_t\in A\}}\,dt.$$
Its Laplace transform is given by
$$\L U(\lambda)=\E\int_0^\infty \exp(-\lambda S_t)\,dt=\dfrac{1}{\varphi(\lambda)}.$$
 By a result of Reveu (see \cite[Proposition 1.7]{Ber2}),
  $U(dx)$ is absolutely continuous with respect to the Lebesgue measure on $[0, \infty)$,
  has a strictly positive bounded continuous density function $u (x)$ on $[0, \infty)$ with $ u(0+)=1$.
 In fact,
\begin{equation}\label{e:3.1a}
  u(x) =   \P (\hbox{there is some $t\geq 0$ so that } S_t=x)  \quad \hbox{for }  x\geq 0.
 \end{equation}
 By \cite[Corollaries 5.4 and 5.5]{BBKR}, $u(x)$ is a completely monotone function on $(0, \infty).$
Let $U^{Z^D}(\cdot, \cdot)$ be the occupation density function of $Z^D.$ That is, $U^{Z^D}(x, y)=\int_0^\infty p^{Z^D}(t, x, y)\,dt,$ where $p^{Z^D}(t, x, y)$ is the transition density function of $Z^D.$ It follows from  \cite[(4.3)]{SV1} that
 \begin{equation}\label{e:3.5'}
 U^{Z^D}(x, y)=\int_0^\infty p^X_D(t, x, y) u(t)dt,
 \end{equation}
 where $p^X_D(t, x, y)$ is the transition density function of the part process $X^D$  killed upon exiting  $D.$

Denote by $G^X_D$  the Green function of $X$ killed upon $D$. Let $W$ be a Brownian motion. Denote by $G^W_D$ the Green function of  $W$ killed upon $D$.
Recall that
$$
g_D(x, y)=\left\{\begin{array}{ll}
\left(1\wedge \dfrac{\delta_D(x) \delta_D(y)}{|x-y|^2}\right)|x-y|^{2-d}, & \quad d\geq 3\\
\log\left(1+\dfrac{\delta_D(x)\delta_D(y)}{|x-y|^2}\right),  & \quad d=2\\
 (\delta_D(x)\delta_D(y))^{1/2} \wedge \dfrac{\delta_D(x)\delta_D(y)}{|x-y|}, & \quad d=1.
\end{array}\right.
$$
When $D$ is a bounded $C^{1, \alpha}$ domain, by \cite[Theorem 4.8]{CKP}, there exists $c=c(d, \lambda_0, \ell, R_0, \Lambda_0, {\rm diam}(D))$ such that
\begin{equation}\label{e:2.5n}
c^{-1}g_D(x, y)\leq  G^X_D(x, y)\leq  cg_D(x, y), \quad x,y\in D.
 \end{equation}
 In particular, when $X$ is a Brownian motion $W,$  there exists $c=c(d,  R_0, \Lambda_0, {\rm diam}(D))$ such that
\begin{equation}\label{e:2.5'}
c^{-1}g_D(x, y)\leq  G^W_D(x, y)\leq  cg_D(x, y), \quad x,y\in D.
 \end{equation}

\begin{prp}\label{P:2.1}
Let $D$ be a bounded  $C^{1, \alpha}$  domain with characteristics $(R_0, \Lambda_0)$ in $\R^d$. There exists a positive constant $C=C(d, \lambda_0, \ell, \phi, R_0, \Lambda_0, {\rm diam}(D))>1$ such that for any $r\in (0, 1),$
$$C^{-1}G^W_{rD}(x, y)\leq U^{Z^{rD}}(x, y)\leq CG^W_{rD}(x, y) \quad {\rm for}\: x,  y\in rD.$$
 \end{prp}

\pf For each $\lambda\geq 1,$ let $X^\lambda_t:=\lambda X_{\lambda^{-2}t}.$
The operator of $X^\lambda_t$ is
$$\L^\lambda f(x)=\nabla (a_{ij}(\lambda^{-1}\cdot) \nabla f)(x).$$
It is easy to see that
 \begin{equation}\label{e:2.5}
 p^X_{\lambda^{-1}D}(t, x, y)=\lambda^{d}p^{X^\lambda}_D(\lambda^2 t, \lambda x, \lambda y), \quad G^X_{\lambda^{-1}D}(x, y)=\lambda^{d-2}G^{X^\lambda}_D(\lambda x, \lambda y), \quad t>0, x, y\in \lambda^{-1}D.
 \end{equation}
 Note that  for each $\lambda\geq 1,$  $a^\lambda_{ij}(\cdot)=a_{ij}(\lambda^{-1}\cdot)$ is $\ell$-Dini continuous and has the uniform elliptic constant $\lambda_0.$
 Thus by \eqref{e:2.5n}, for any $\lambda\geq 1,$
 $G^{X^\lambda}_D(x, y)\asymp G^W_D(x, y)\asymp g_D(x, y),$ where the comparison constants depend on $(d, \lambda_0, \ell, R_0, \Lambda_0, {\rm diam}(D)).$
Hence, by this comparability and \eqref{e:2.5},
  \begin{equation}\label{e:2.6}
  G^X_{\lambda^{-1}D}(x, y)\asymp G^W_{\lambda^{-1}D}(x, y)\asymp g_{\lambda^{-1}D}(x, y), \quad x, y\in \lambda^{-1}D,
  \end{equation}
  where the comparison constants depend only on $(d, \lambda_0, \ell, R_0, \Lambda_0, {\rm diam}(D)).$

Since the potential density function $u$ of $S_t$ is
bounded by $1$,
we have by  \eqref{e:3.5'} and \eqref{e:2.6}, there exists $c_1=c_1(d, \lambda_0, \ell, R_0, \Lambda_0, {\rm diam}(D))$ such that for any $r\in (0, 1),$
\begin{equation}\label{e:2.8}
 U^{Z^{rD}}(x, y)\leq \int_0^\infty p^X_{rD}(t, x, y) dt= G^X_{rD}(x, y)\leq c_1 G^W_{rD}(x, y), \quad {\rm for}\: x,  y\in rD.
 \end{equation}
Fix $r\in (0, 1).$ Let $\lambda:=r^{-1}.$ Note that    $a^\lambda_{ij}(\cdot)=a_{ij}(\lambda^{-1}\cdot)$ is $\ell$-Dini continuous and has the uniform elliptic constant $\lambda_0.$ By the result in \cite{Cho}, for fixed $T>0,$ there exist positive constants $c_k=c_k(d, \lambda_0, \ell, R_0, \Lambda_0,  T), k=2, 3$ such that for any $ x, y\in D$ and $t\in (0, T),$
$$
p^{X^\lambda}_D(t, x, y)\geq c_2 \left(1\wedge \dfrac{\delta_D(x) \delta_D(y)}{t}\right)t^{-d/2} \exp(-c_3\dfrac{|x-y|^2}{t}).
$$
By this inequality and \eqref{e:2.5}, for each $T>0,$
 \begin{equation}\label{e:2.10n}
p^{X}_{rD}(t, x, y)=r^{-d}p^{X^\lambda}_D(r^{-2} t, r^{-1} x, r^{-1} y)\geq c_2 \left(1\wedge \dfrac{\delta_{rD}(x) \delta_{rD}(y)}{t}\right)t^{-d/2} \exp(-c_3\dfrac{|x-y|^2}{t})
 \end{equation}
 for any $x, y\in rD, t\leq r^2T.$
By taking $T={\rm diam}(D)^2$ in \eqref{e:2.10n}, we have
\begin{eqnarray}\label{e:2.9n}
&&\int_0^{( {\rm diam} (rD))^2} p^X_{rD}(t, x, y) dt \nonumber\\
&\geq & c_2 \int_0^{({\rm diam} (rD))^2}\left(1\wedge  \dfrac{\delta_{rD}(x)\delta_{rD}(y)}{t}\right) t^{-d/2} \exp(-c_3\dfrac{|x-y|^2}{t})\,dt.\label{e:2.12'}
\end{eqnarray}

When $d\geq 3,$ by \eqref{e:2.9n} and the change of variable $s=|x-y|^2/t$,
\begin{eqnarray}\label{e:2.10}
&&\int_0^{( {\rm diam} (rD))^2} p^X_{rD}(t, x, y) dt \nonumber\\
&\geq &\int_0^{|x-y|^2} p^X_{rD}(t, x, y) dt \nonumber\\
&\geq & c_4|x-y|^{2-d}\left(1\wedge \dfrac{\delta_{rD}(x)\delta_{rD}(y)}{|x-y|^2}\right)\int_{1}^\infty  s^{d/2-2} \exp(-c_3s) \,ds\nonumber \\
&\geq & c_5|x-y|^{2-d}\left(1\wedge \dfrac{\delta_{rD}(x)\delta_{rD}(y)}{|x-y|^2}\right)
= c_5g_{rD}(x, y).
\end{eqnarray}

 When $d=1, 2,$  by \eqref{e:2.9n} and a very similar argument in \cite[Proposition 3.3 and Remark 3.4]{KSV1} and \cite[Theorem 4.1]{CKSV1} respectively,  there exists $c_6$ such that for any $r\in (0, 1),$
\begin{equation}\label{e:2.11}
\int_0^{( {\rm diam} (rD))^2} p^X_{rD}(t, x, y) dt \geq c_6 g_{rD}(x, y).
\end{equation}

Since $u(t)$ is  continuous decreasing and strictly positive with $ u(0+)=1$,
then $u(t)\geq u({\rm diam}(D)^2)=:c_7$ for $t\in (0, ({\rm diam}(D))^2).$ By combining \eqref{e:2.10}-\eqref{e:2.11},
$$
 U^{Z^{rD}}(x, y) =\int_0^\infty p^X_{r D}(t, x, y) u(t)dt
 \geq c_7\int_0^{({\rm diam} (rD))^2} p^X_{r D}(t, x, y) dt\geq c_8g_{rD}(x, y)
  \geq c_9G^W_{rD}(x, y),
$$
where the last inequality is due to \eqref{e:2.6}.
This together with \eqref{e:2.8}   proves the result.
 \qed

\medskip
Let $\psi(r):=H(r^{-2})^{-1}.$ If $H$ satisfies $L^a(\gamma, c_L)$ and $U^a(\delta, C_U)$ for some $a>0$ with  $\delta<2,$
then $\psi$ is a non-negative function  satisfying $L_{a^{-1/2}}(2\gamma, c_L)$ and $U_{a^{-1/2}}(2\delta, C_U).$ It follows from \eqref{e:2.16} that when $H$ satisfies $L^a(\gamma, c_L)$ and $U^a(\delta, C_U)$ for some $a>0$ with  $\delta< 2,$ for each $M>0,$ there exists a constant $c$ depending on $M$ such that
 \begin{equation}\label{e:2.9}
 c^{-1}\dfrac{1}{|x-y|^d\psi(|x-y|)} \leq j(|x-y|)\leq c\dfrac{1}{|x-y|^d\psi(|x-y|)}, \quad x, y\in\R^d\setminus {\rm diag}\quad {\rm with} \quad |x-y|\leq M.
 \end{equation}
Note that $\int_{\R^d} (1\wedge |z|^2) j(|z|)\,dz<\infty.$
Hence by \eqref{e:2.9},
\begin{equation}\label{e:1.12}
\int_0^1 s/\psi(s)\,ds<\infty.
\end{equation}
Let $\psi_0(r):=\frac{r^2}{\psi(r)}.$ \eqref{e:2.9} is equivalent that for each $M>0,$ there exists a constant $c$ depending on $M$ such that
\begin{equation}\label{e:2.9'}
 c^{-1}\dfrac{\psi_0(|x-y|)}{|x-y|^{d+2}} \leq j(|x-y|)\leq c\dfrac{\psi_0(|x-y|)}{|x-y|^{d+2}}, \quad x, y\in\R^d\setminus {\rm diag}\quad {\rm with} \quad |x-y|\leq M.
 \end{equation}

\begin{lem}\label{L:2.3n}
Suppose $H$ satisfies $L^a(\gamma, c_L)$ and $U^a(\delta, C_U)$ for some $a>0$ with  $\delta<2.$
 Then  there exists $c=c(a)>0$ such that for any  $r\in (0, a^{-1/2}),$
\begin{equation}\label{e:2.23}
c^{-1}\dfrac{1}{\psi(r)}\leq \int_r^{a^{-1/2}} \dfrac{1}{s\psi(s)}\,ds\leq c\dfrac{1}{\psi(r)}.
\end{equation}
\end{lem}

\pf If $H$ satisfies $L^a(\gamma, c_L)$ and $U^a(\delta, C_U)$ for some $a>0$ with  $\delta<2$, then for each $M>0,$ $\psi$ is a non-negative function  satisfying $L_{a^{-1/2}}(2\gamma, c_L)$ and $U_{a^{-1/2}}(2\delta, C_U).$
Note that
$$\int_r^{a^{-1/2}} \dfrac{1}{s\psi(s)}\,ds= \dfrac{1}{\psi(r)}\int_r^{a^{-1/2}} \dfrac{\psi(r)}{s\psi(s)}\,ds.$$
Since  $\psi$ is a non-negative function  satisfying $L_{a^{-1/2}}(2\gamma, c_L)$ and $U_{a^{-1/2}}(2\delta, C_U),$ then it is easy to obtain \eqref{e:2.23}.
\qed

\begin{lem}\label{L:2.2}
 If $H$ satisfies $L^a(\gamma, c_L)$ and $U^a(\delta, C_U)$ for some $a>0$ with  $\delta<2$, then
 \begin{equation}\label{e:2.4}
c_L(\dfrac{r}{R})^{2-2\gamma}\leq \dfrac{\psi_0(r)}{\psi_0(R)}\leq C_U(\dfrac{r}{R})^{2-2\delta} \quad \mbox{for} \quad 0<r<R\leq a^{-1/2}
\end{equation}
and
\begin{equation}\label{e:2.23'}
\lim_{s\rightarrow 0} \psi_0(s)=0.
\end{equation}
\end{lem}

\pf If $H$ satisfies $L^a(\gamma, c_L)$ and $U^a(\delta, C_U)$ with  $\delta<2,$ then $\psi$ is a non-negative function  satisfying $L_{a^{-1/2}}(2\gamma, c_L)$ and $U_{a^{-1/2}}(2\delta, C_U).$
Note that $\psi_0(r)/\psi_0(R)= \frac{\psi(R)}{\psi(r)}(\frac{r}{R})^2.$
Then it is easy to see that \eqref{e:2.4} holds.

By \eqref{e:2.4},  for $s\in (0, a^{-1/2}/2),$
\begin{equation}\label{e:2.20'}
C_U^{-1}2^{-(3-2\delta)}\psi_0(s)\leq \int_s^{2s} \dfrac{\psi_0(r)}{r}\,dr\leq c_L^{-1}2^{3-2\gamma}\psi_0(s).
\end{equation}
It follows from \eqref{e:1.12} that $\int_0^1 s/\psi(s)\,ds<\infty.$ Hence $\int_0^1 \frac{\psi_0(r)}{r}\,dr<\infty.$
Thus
$$\lim_{s\rightarrow 0}\int_s^{2s} \dfrac{\psi_0(r)}{r}\,dr=0.$$
Consequently, \eqref{e:2.23'} holds by \eqref{e:2.20'}.

\qed

\smallskip

The following  resurrection formula is from the combination of Theorem 4.1 and Corollaries 4.2-4.3 in \cite{SV1}.

\begin{thm}\label{T:2.2}
Let $D$ be an open set.
For each open set $B\subset \overline D$ and $C\subset \R^d,$
$$\begin{aligned}
&\P_x(Y_{\sigma_{\tau^X_D}-}\in B, Y_{\sigma_{\tau^X_D}}\in C)\\
=& \int_{B\cap D} U^{Z^D}(x, y)\int_{C\cap D^c} J(y, z)\,dz\,dy
+\int_{B\cap D} U^{Z^D}(x, y)\int_{C\cap D} (J(y, z)-J^{Z^D}(y, z))\,dz\,dy\\
&+\E_x [u(\tau^X_D); X_{\tau^X_D-}\in B, X_{\tau^X_D}\in C\cap \partial D], \quad x\in D.
\end{aligned}$$
In particular, for each Borel set $C\subset \partial D,$
$$\P_x(Y_{\sigma_{\tau^X_D}}\in C; S_{\sigma_{\tau^X_D}}=\tau^X_D)=\E_x [u(\tau^X_D); X_{\tau^X_D}\in C], \quad x\in D.$$
\end{thm}

\smallskip

Define
\begin{equation}\label{e:2.21'}
q_{D}(y, z):=J(y, z)-J^{Z^{D}}(y, z), \quad y, z\in D.
\end{equation}
By \eqref{e:1.6n} and \eqref{e:2.2},
\begin{equation}\label{e:2.10'}
q_{D}(y, z)=\int_0^\infty (p^X(t, y, z)-p^X_{D}(t, y, z))\,\mu(t)dt,
\end{equation}
Denote by $p^X(t, x, y)$ and $p^X_{D}(t, x, y)$  the transition density functions of diffusion $X$ in $\R^d$ and the subprocess $X^D$ in $D.$
In the remainder of this paper, we always use the constant $a$ to denote the constant in the assumption (A1).

\begin{lem}\label{L:2.4}
Suppose $H$ satisfies the assumption (A1).
 Suppose $D$ is a bounded $C^{1, \alpha}$ domain with characteristics $(R_0, \Lambda_0),$ there exists $C=C(d, \lambda_0, \ell, \phi,  R_0, \Lambda_0, {\rm diam}(D))>1$  such that for any $r\in (0, \frac{a^{-1/2}}{2{\rm diam}(D)}),$
$$\int_{rD}\int_{rD} G^W_{rD}(x, y) q_{rD}(y, z)\delta_{rD}(z)\,dy\,dz\leq C\Psi(r{\rm diam}(D))\cdot \delta_{rD}(x), \quad x\in D,$$
where $\Psi(r)\rightarrow 0$ as $r\rightarrow 0.$
\end{lem}

\pf Without loss of generality, we assume the constant  $a=1$ in the assumption (A1). Let $r\in (0, \frac{1}{2{\rm diam}(D)}).$ By \eqref{e:2.10'} and \eqref{e:2.8'}, there exists $c_1=c_1(d, \lambda_0)$ such that
\begin{equation}\label{e:2.11'}
q_{rD}(y, z)\leq \int_0^\infty p^X(t, y, z)\,\mu(t)dt=  J(y, z)\leq c_1j(|y-z|), \quad y, z\in rD.
\end{equation}
Since
 $$p^X(t, y, z)-p^X_{rD}(t, y, z)=\E_y [p^X(t-\tau^X_{rD}, X_{\tau^X_{rD}}, z); \tau^X_{rD}<t]$$
 and $\mu(t)$ is a decreasing function,
we have for $y, z\in D,$
\begin{equation}\label{e:2.12}
\begin{aligned}
 q_{rD}(y, z)&=\int_0^\infty \E_y [p^X(t-\tau^X_{rD}, X_{\tau^X_{rD}}, z); \tau^X_{rD}<t]\,\mu(t)dt\\
 &=\int_{(0, \infty)\times \partial D} \int_s^\infty  p^X(t-s, u, z)\mu(t)\,dt \cdot\P_y((\tau^X_{rD}, X_{\tau^X_{rD}})\in (ds, du))\\
 &=\int_{(0, \infty)\times \partial D} \int_0^\infty  p^X(t, u, z)\mu(t+s) dt \cdot\P_y((\tau^X_{rD}, X_{\tau^X_{rD}})\in (ds, du))\\
 &\leq \int_{(0, \infty)\times \partial D} \int_0^\infty  p^X(t, u, z)\mu(t)\,dt \cdot\P_y((\tau^X_{rD}, X_{\tau^X_{rD}})\in (ds, du))\\
   &=\int_{(0, \infty)\times \partial D} J(u,z) \cdot\P_y((\tau^X_{rD}, X_{\tau^X_{rD}})\in (ds, du))\\
 &= \E_y J(X_{\tau^X_{rD}}, z)\\
 &\leq c_2 j(\delta_{rD}(z)),
 \end{aligned}\end{equation}
 where the last inequality is due to \eqref{e:2.8'} and that $j(r)$ is a decreasing function.
 Hence, by \eqref{e:2.11'} and \eqref{e:2.12}, there exists $c_3=c_3(d, \lambda_0)$ such that
 \begin{equation}\label{e:2.13}
 q_{rD}(y, z)\leq c_3\left(j(|y-z|)\wedge j(\delta_{rD}(z))\right), \quad y, z\in rD.
 \end{equation}
It follows from \eqref{e:2.16} that there exists $c_4$ such that $j(s)\leq c_4j(2s)$ for $s\in (0, 1/2).$
 Note that ${\rm diam}(rD)\leq \frac{1}{2}.$
By \eqref{e:2.13}, \eqref{e:2.9} and Lemma \ref{L:2.3n}, we have
\begin{equation}\label{e:2.16'}\begin{aligned}
&\int_{rD} q_{rD}(y, z)\delta_{rD}(z)\,dz\\
\leq &c_3\int_{rD\cap B(y, \delta_{rD}(y)/2)} j(\delta_{rD}(z))\delta_{rD}(z)\,dz
+c_3\int_{rD\setminus B(y, \delta_{rD}(y)/2)} j(|y-z|) \delta_{rD}(z)\,dz\\
\leq &c_5\int_{rD\cap B(y, \delta_{rD}(y)/2)}\delta_{rD}(y)j(\delta_{rD}(y)) \,dz\\
&+2c_3\delta_{rD}(y)\int_{rD\cap \{z: |y-z|>\delta_{rD}(y)/2, \delta_{rD}(z) \leq 2\delta_{rD}(y)\}} j(|y-z|)\,dz\\
&+3c_3\int_{rD\cap \{z: |y-z|>\delta_{rD}(y)/2, \delta_{rD}(z) > 2\delta_{rD}(y)\}} j(|y-z|) |y-z|\,dz\\
\leq & c_6\delta_{rD}(y)j(\delta_{rD}(y))(\delta_{rD}(y))^d+c_6\delta_{rD}(y)\int_{rD\cap \{z: |y-z|>\delta_{rD}(y)/2, \delta_{rD}(z) \leq 2\delta_{rD}(y)\}} \dfrac{1}{|y-z|^d\psi(|y-z|)}\,dz\\
&+c_6\int_{rD\cap \{z: |y-z|>\delta_{rD}(y)/2\}}\dfrac{1}{|y-z|^{d-1}\psi(|y-z|)}\,dz\\
\leq & c_6\dfrac{\delta_{rD}(y)}{\psi(\delta_{rD}(y))}+c_6\int_{\{\delta_{rD}(y)/2<s<{\rm diam}(rD)\}} \left(\dfrac{\delta_{rD}(y)}{s\psi(s)}+\dfrac{1}{\psi(s)}\right)\,ds\\\leq & c_7\dfrac{\delta_{rD}(y)}{\psi(\delta_{rD}(y))}+c_7\int_{\{\delta_{rD}(y)/2<s<{\rm diam}(rD)\}} \dfrac{1}{\psi(s)}\,ds
\end{aligned}\end{equation}
where in the second  inequality, we used  $\delta_{rD}(z)\leq |z-y|+\delta_{rD}(y)\leq 3|z-y|$ for $z\in rD$ with $|z-y|>\delta_{rD}(y)/2.$

In the following, we estimate
\begin{equation}\label{e:2.26n}\begin{aligned}
&\int_{rD}\int_{rD} G^W_{rD}(x, y) q_{rD}(y, z)\delta_{rD}(z)\,dy\,dz\\
=&\int_{rD\cap B(x, \delta_{rD}(x)/2)}\int_{rD} G^W_{rD}(x, y) q_{rD}(y, z)\delta_{rD}(z)\,dz\,dy\\
&+\int_{rD\cap B^c(x, \delta_{rD}(x)/2)}\int_{rD} G^W_{rD}(x, y) q_{rD}(y, z)\delta_{rD}(z)\,dz\,dy\\
=:&I+II.
\end{aligned}\end{equation}

For the first term $I,$
recall that $\psi_0(s)= s^2 /\psi(s).$
By \eqref{e:2.4},
\begin{equation}\label{e:2.17'}
\begin{aligned}
&\int_{\{\delta_{rD}(y)/2<s<{\rm diam}(rD)\}} \dfrac{1}{\psi(s)}\,ds
=\int_{\{\delta_{rD}(y)/2<s<{\rm diam}(rD)\}} \dfrac{\psi_0(s)}{s^2}\,ds\\
&\leq C_U\psi_0({\rm diam}(rD))\int_{\delta_{rD}(y)/2}^{{\rm diam}(rD)}\dfrac{1}{s^2}\,ds
\leq 2C_U\delta_{rD}(y)^{-1} \psi_0({\rm diam}(rD)).
\end{aligned}\end{equation}
 By \eqref{e:2.16'} and \eqref{e:2.17'}, we have
\begin{equation}\label{e:2.30'}\begin{aligned}
I&\leq (c_7+2C_U)\int_{rD\cap B(x, \delta_{rD}(x)/2)}G^W_{rD}(x, y)
\left(\dfrac{\delta_{rD}(y)}{\psi(\delta_{rD}(y))}+\delta_{rD}(y)^{-1}\psi_0({\rm diam}(rD))\right)\,dy\\
&\leq c_8\left(\dfrac{\delta_{rD}(x)}{\psi(\delta_{rD}(x))}+\delta_{rD}(x)^{-1}\psi_0({\rm diam}(rD))\right)\int_{rD\cap B(x, \delta_{rD}(x)/2)} G^W_{rD}(x, y) \,dy.
\end{aligned}\end{equation}
It follows from \eqref{e:2.6} that $G^W_{rD}(x, y)\leq c_9g_{rD}(x, y).$
Then it is easy to calculate that for $d\geq 1,$
\begin{equation}\label{e:2.31'}
\int_{B(x, \delta_{rD}(x)/2)}G^W_{rD}(x, y) \,dy\leq c_9 \int_{B(x, \delta_{rD}(x)/2)}g_{rD}(x, y) \,dy\leq c_{10}(\delta_{rD}(x))^2.
\end{equation}
Hence, by \eqref{e:2.30'} and \eqref{e:2.31'},
\begin{equation}\label{e:2.32n}\begin{aligned}
I&\leq  c_8c_{10}\delta_{rD}(x)\left[\frac{\delta^2_{rD}(x)}{\psi(\delta_{rD}(x))}+\psi_0({\rm diam}(rD))\right]\\
&=c_8c_{10}\delta_{rD}(x)(\psi_0(\delta_{rD}(x))+\psi_0({\rm diam}(rD)))\\
&\leq c_{11}\delta_{rD}(x)\psi_0({\rm diam}(rD)),
\end{aligned}\end{equation}
where in the last inequality we used \eqref{e:2.4}.

Next we estimate the second term $II$ in \eqref{e:2.26n}.
We assert that there exist $c_{12}>1$ and  $\ee\in (1\vee(2\delta), 2)$  such that
\begin{equation}\label{e:3.27n}
\int_s^1 \dfrac{1}{\psi(u)}\,du\leq c_{12}(\frac{\psi_0(s)}{s}+s^{1-\ee}) \quad \mbox{for} \quad s\in (0, 1].
\end{equation}
In fact, if $H$ satisfies $U^a(\delta, C_U)$ with $\delta<1,$ then  $\psi$ is a non-negative function  satisfying $U_{a^{-1/2}}(2\delta, C_U).$
Hence, $\psi(s)^{-1}\leq c_{13}s^{-2\delta}$ for $s\in (0, 1).$
Let $\psi_1(s):=\int_s^1 \frac{1}{\psi(u)}\,du.$
We have
 $$
\psi_1(s)\leq c_{13}\int_s^1 u^{-2\delta}\,du \leq c_{14}(s^{1-2\delta}1_{2\delta\in (1, 2)}+\log(s^{-1})1_{2\delta\leq 1}) \quad \mbox{for} \quad s\in (0, 1).
$$
Observe that for any  $\ee\in (1, 2)$, $\log(s^{-1})\leq s^{1-\ee}$ for $s\in (0, 1).$
Thus there exists $\ee\in (1\vee(2\delta), 2)$ such that
\begin{equation}\label{e:2.35'}
\psi_1(s)\leq  c_{14}s^{1-\ee}\quad \mbox{for} \quad s\in (0, 1).
\end{equation}
On the other hand, if  $H$ satisfies $L^a(\gamma, c_L)$ and $U^a(\delta, c_L)$ with $\delta=1$ and $\gamma>1/2,$ then $\psi$  satisfies $L_{a^{-1/2}}(2\gamma, c_L)$ with $\gamma>1/2.$
Then
\begin{equation}\label{e:2.22'}
\int_s^1 \dfrac{1}{\psi(u)}\,du=\dfrac{s}{\psi(s)}\int_s^1 s^{-1}\dfrac{\psi(s)}{\psi(u)}\,du\leq \dfrac{s}{\psi(s)}\int_s^1 c_L^{-1}s^{-1}\dfrac{s^{2\gamma}}{u^{2\gamma}}\,du\leq \frac{c_L^{-1}}{2\gamma-1} \dfrac{s}{\psi(s)}=\frac{c_L^{-1}}{2\gamma-1}\dfrac{\psi_0(s)}{s}.
\end{equation}
Hence, \eqref{e:3.27n} follows by \eqref{e:2.35'} and \eqref{e:2.22'}.

Let $\ee\in (1\vee(2\delta), 2)$ be the constant in \eqref{e:3.27n}.
We divide three cases to estimate the second term $II$ in \eqref{e:2.26n}.
When $d\geq 3,$
by \eqref{e:2.6}, \eqref{e:2.16'} and \eqref{e:3.27n},
\begin{equation}\label{e:2.34n}\begin{aligned}
II&\leq c_{15}\int_{rD\cap B^c(x, \delta_{rD}(x)/2)}G^W_{rD}(x, y)\left(\dfrac{\psi_0(\delta_{rD}(y))}{\delta_{rD}(y)}+
(\delta_{rD}(y))^{1-\ee}\right)\,dy\\
&\leq c_{16}\int_{rD\cap B^c(x, \delta_{rD}(x)/2)\}}\dfrac{\delta_{rD}(x)}{|x-y|}|x-y|^{2-d}\dfrac{\delta_{rD}(y)}{|x-y|}
\left(\dfrac{\psi_0(|x-y|)}{\delta_{rD}(y)}+(\delta_{rD}(y))^{1-\ee} \right)\,dy\\
&\leq c_{17}\delta_{rD}(x)\int_{rD\cap B^c(x, \delta_{rD}(x)/2)\}}\left( \dfrac{\psi_0(|x-y|)}{|x-y|^d}+|x-y|^{-d}(\delta_{rD}(y))^{2-\ee}\right)\,dy\\
&\leq c_{18}\delta_{rD}(x)\int_{rD\cap B^c(x, \delta_{rD}(x)/2)\}}\left( \dfrac{\psi_0(|x-y|)}{|x-y|^d}+|x-y|^{2-\ee-d}\right)\,dy\\
&\leq c_{19} \delta_{rD}(x)\left[\int_0^{{\rm diam}(rD)} \dfrac{\psi_0(s)}{s}\,ds+{\rm diam}(rD)^{2-\ee}\right],\\
\end{aligned}\end{equation}
where in the second and fourth inequalities, we used  $\delta_{rD}(y)\leq |y-x|+\delta_{rD}(x)\leq 3|y-x|$ for $y\in D\cap B^c(x, \delta_{rD}(x)/2)$ and \eqref{e:2.4}.

When $d=2,$ we have by \eqref{e:2.6}, \eqref{e:2.16'} and \eqref{e:3.27n},
\begin{equation}\label{e:2.35n}\begin{aligned}
II&\leq c_{20}\int_{rD\cap \{y: |y-x|> 2\delta_{rD}(x)\}}G^W_{rD}(x, y)\left(\delta_{rD}(y)^{-1}\psi_0(\delta_{rD}(y))+(\delta_{rD}(y))^{1-\ee}\right)dy\\
&\leq c_{21}\int_{rD\cap \{y: |y-x|> 2\delta_{rD}(x)\}}\log(1+\dfrac{\delta_{rD}(x)\delta_{rD}(y)}{|x-y|^2})
\left(\delta_{rD}(y)^{-1}\psi_0(\delta_{rD}(y))+(\delta_{rD}(y))^{1-\ee}\right)\,dy\\
&\leq c_{22}\int_{rD\cap \{y: |y-x|> 2\delta_{rD}(x)\}} \left(\dfrac{\delta_{rD}(x)}{|x-y|}\dfrac{\delta_{rD}(y)}{|x-y|}
\delta_{rD}(y)^{-1}\psi_0(\delta_{rD}(y))+\dfrac{\delta_{rD}(x)}{|x-y|^2}(\delta_{rD}(y))^{2-\ee}\right)\,dy\\
&\leq c_{23}\delta_{rD}(x)\int_{rD\cap \{y: |y-x|> 2\delta_{rD}(x)\}} (|x-y|^{-2}\psi_0(|x-y|)+|x-y|^{-\ee})\,dy\\
&\leq c_{24}\delta_{rD}(x)\left[\int_0^{{\rm diam}(rD)} \dfrac{\psi_0(s)}{s}\,ds+{\rm diam}(rD)^{2-\ee}\right].
\end{aligned}\end{equation}

When $d=1,$ by \eqref{e:2.6}, \eqref{e:2.16'} and \eqref{e:3.27n},
\begin{equation}\label{e:2.36n}\begin{aligned}
II&\leq c_{25}\int_{rD\cap \{y: |y-x|> 2\delta_{rD}(x)\}}G^W_{rD}(x, y)\left(\delta_{rD}(y)^{-1}\psi_0(\delta_{rD}(y))+(\delta_{rD}(y))^{1-\ee}\right)dy\\
&\leq c_{26}\delta_{rD}(x)\int_{D\cap B^c(x, \delta_{rD}(x)/2)}\dfrac{\delta_{rD}(y)}{|x-y|}\left(\delta_{rD}(y)^{-1}\psi_0(\delta_{rD}(y))+(\delta_{rD}(y))^{1-\ee}\right)\,dy\\
&\leq c_{27}\delta_{rD}(x)\int_{rD\cap B^c(x, \delta_{rD}(x)/2)}(|x-y|^{-1}\psi_0(|x-y|)+|x-y|^{1-\ee})\,dy\\
&\leq c_{28}\delta_{rD}(x)\left[\int_0^{{\rm diam}(rD)} \dfrac{\psi_0(s)}{s}\,ds+{\rm diam}(rD)^{2-\ee}\right].
\end{aligned}\end{equation}

Define $\Psi(r):=\psi_0(r)+\int_0^r \psi_0(s)/s\,ds+r^{2-\ee}.$
By combing \eqref{e:2.26n}, \eqref{e:2.32n}, \eqref{e:2.34n}-\eqref{e:2.36n}, there exists $c_{29}=c_{29}(d, \lambda_0, \ell, \phi,  R_0, \Lambda_0, {\rm diam}(D))>1$ such that
$$\int_{rD}\int_{rD} G^W_{rD}(x, y) q_{rD}(y, z)\delta_{rD}(z)\,dy\,dz\leq c_{29}\delta_{rD}(x)\Psi({\rm diam}(rD)).$$
Due to $\int_0^1 \frac{\psi_0(s)}{s}\,dr<\infty$ and \eqref{e:2.23'},  $\Psi(r)\rightarrow 0$ as $r\rightarrow 0.$ Hence, the desired conclusion is obtained.

\qed

\smallskip

Recall that  $p(t, x, y)$ is the transition density function of $Y$ in $\R^d.$ For each open set $B,$ denote by $\tau_B$ the first exiting time of $Y$ from $B.$

\begin{lem}\label{L:2.5'} 
Suppose $H$ satisfies $L^a(\gamma, c_L)$ and $U^a(\delta, C_U)$ with  $\delta<2$ for some $a>0.$
There exists $C=C(d, \lambda_0, \ell, \phi)$ such that for any $x_0\in\R^d$ and $r\in (0, 1),$
$$C^{-1}r^2\leq \inf_{x\in B(x_0, r/2)}\E_x \tau_{B(x_0, r)}\leq \sup_{x\in B(x_0, r)}\E_x \tau_{B(x_0, r)}\leq C r^2.$$
\end{lem}

\pf Fix $x_0\in\R^d.$ For the simplicity of notation, let $B_r:=B(x_0, r).$
The proof of upper bound is standard (see e.g. \cite[Lemma 2.3]{CK3}).
By Theorem \ref{T0},
$$p(t, x, y)\leq c_1t^{-d/2} \quad \mbox{for}\quad t\in (0, 1), x, y\in B_r.$$
We choose $c_2>0$ such that $c_1(c_2r^2)^{-d/2}m_d(B_r)\leq 1/2.$ Let $t:=c_2r^2.$ Then
$$\P_x(Y_t\in B(x_0, r))=\int_{B_r} p(t, x, y)\,dy\leq 1/2, \quad x\in B_r.$$
Hence for each $x\in B_r,$ $\P_x(\tau_{B_r}\leq t)\geq \P_x(Y_t\in B^c_r) \geq 1/2.$ That is $\P_x(\tau_{B_r}> t)\leq 1/2.$
Then by the strong Markov property of $Y$ and the induction argument,
$\P_x(\tau_{B_r}>kt)\leq 2^{-k}$ for each $k\geq 1.$
This yields that $\sup_{x\in B_r}\E_x \tau_{B_r}\leq \sup_{x\in B_r}\sum_{k=0}^\infty t\P_x(\tau_{B_r}> kt)\leq c_3r^2.$

For the lower bound,
let $Z^{B_r}_t:=X^{B_r}_{S_t}$ be  the subordinate killed diffusion in $B_r$. Let $\zeta$  denote the life time of the process $Z^{B_r}_t.$
Then by Proposition \ref{P:2.1}, for $x\in B_{r/2},$
$$\E_x \tau_{B_r}\geq \E_x \zeta=\int_{B_r} U^{Z^{B_r}}(x, y)\,dy\geq c_4\int_{B_r} G^W_{B_r}(x, y)\,dy=c_4\E_x \tau^W_{B_r}
\geq c_5\E_x \tau^W_{B(x, r/4)}\geq c_6r^2.$$
Hence, the proof is complete.

\qed

\begin{prp}\label{L:2.5}
Suppose $H$ satisfies the assumption (A1).
 Suppose $D$ is a bounded $C^{1, \alpha}$ domain with characteristics $(R_0, \Lambda_0),$ there exist positive constants $\delta_1=\delta_1(d, \lambda_0, \ell,  \phi, R_0, \Lambda_0, {\rm diam}(D))\in (0, R_0)$ and
 $C=C(d, \lambda_0, \ell,  \phi, R_0, \Lambda_0, {\rm diam}(D))$ such that for any $r\in (0, \delta_1/{\rm diam}(D)),$
$$C^{-1} \E_x\tau^W_{rD}\leq \E_x\tau_{rD}\leq C\E_x\tau^W_{rD}, \quad x\in rD.$$
\end{prp}

\pf Let  $Z^{rD}_t:=X^{rD}_{S_t}$ be  the subordinate killed diffusion in $rD$. Let $\zeta$  denote the life time of the process $Z^{rD}_t.$
It follows from  Proposition \ref{P:2.1} that there exists $c_1=c_1(d, \lambda_0, \ell, \phi, R_0, \Lambda_0, {\rm diam}(D))>0$ such that for any $r\in (0, \frac{a^{-1/2}}{2{\rm diam}(D)}),$
\begin{equation}\label{e:2.19'}
\E_x \tau_{rD}\geq \E_x \zeta=\int_{rD} U^{Z_{rD}}(x, y)\,dy\geq c_1\int_{rD} G^W_{rD}(x, y)\,dy= c_1\E_x\tau^W_{rD}.
\end{equation}
Let  $\tau^X_{rD}$ be the first exiting time of $X$ from $B.$
Note that $\zeta=\sigma_{\tau^X_{rD}}=\inf\{t>0: S_t>\tau^X_{rD}\}.$
By the strong Markov property of $Y,$ we have
\begin{equation}\label{e:2.19}\begin{aligned}
\E_x \tau_{rD}&= \int_0^\infty \P_x(\tau_{rD}>t)\,dt\\
&=\int_0^\infty \P_x(\sigma_{\tau^X_{rD}}>t)\,dt+ \int_0^\infty \P_x(\tau_{rD}>t\geq \sigma_{\tau^X_{rD}})\,dt\\
&=\int_0^\infty \P_x(\zeta>t)\,dt+\int_0^\infty \P_x(Y_{\sigma_{\tau^X_{rD}-}}\in rD; Y_{\sigma_{\tau^X_{rD}}}\in rD; \tau_{rD}>t\geq\sigma_{\tau^X_{rD}})\,dt\\
&= \E_x\zeta+
\int_0^\infty \E_x(Y_{\sigma_{\tau^X_{rD}-}}\in rD, Y_{\sigma_{\tau^X_{rD}}}\in rD,  t\geq\sigma_{\tau^X_{rD}}; \P_{ Y_{\sigma_{\tau^X_{rD}}}}(\tau_{rD}>t-\sigma_{\tau^X_{rD}}))\,dt\\
&=\E_x\zeta+
\E_x\left[Y_{\sigma_{\tau^X_{rD}-}}\in D; Y_{\sigma_{\tau^X_{rD}}}\in rD; \E_{ Y_{\sigma_{\tau^X_{rD}}}}\tau_{rD}\right], \quad x\in rD.
\end{aligned}\end{equation}
By Theorem \ref{T:2.2}, we have
\begin{equation}\label{e:2.20}\begin{aligned}
&\E_x(Y_{\sigma_{\tau^X_{rD}-}}\in rD; Y_{\sigma_{\tau^X_{rD}}}\in rD; \E_{ Y_{\sigma_{\tau^X_{rD}}}}\tau_{rD})\\
=&\int_{rD}\int_{rD} U^{Z^{rD}}(x, y) q_{rD}(y, z) \E_z \tau_{rD}\,dy\,dz\\
\end{aligned}\end{equation}
Thus by \eqref{e:2.19}-\eqref{e:2.20}, for $x\in rD,$
\begin{equation}\label{e:2.41}
\E_x \tau_{rD}=\E_x\zeta+\int_{rD}\int_{rD} U^{Z^{rD}}(x, y) q_{rD}(y, z) \E_z \tau_{rD}\,dy\,dz.
\end{equation}
For the simplicity of notation, denote by $(U^{Z^{rD}}\ast  q_{rD}) (x, z)=\int_{rD} U^{Z^{rD}}(x, y) q_{rD}(y, z)\,dy.$
By \eqref{e:2.41} and the induction, we have
\begin{equation}\label{e:2.21}\begin{aligned}
\E_x \tau_{rD}&=\E_x\zeta + \sum_{n=1}^N\int_{rD} (U^{Z^{rD}}\ast  q_{rD})^n (x, z) \E_z \zeta\,dz\\
&+ \int_{rD} (U^{Z^{rD}}\ast  q_{rD})^{N+1} (x, z) \E_z \tau_{rD}\,dz.\\
\end{aligned}\end{equation}
It follows from Theorem \ref{T:2.2} that
\begin{equation}\label{e:2.47}\begin{aligned}
\int_{rD} U^{Z^{rD}}\ast  q_{rD} (x, z)\,dz
&=\P_x(Y_{\sigma_{\tau^X_{rD}}-}\in rD, Y_{\sigma_{\tau^X_{rD}}}\in rD)\\
&\leq 1-\P_x(\tau^X_{rD}=S_{\sigma_{\tau^X_{rD}}}; Y_{\sigma_{\tau^X_{rD}}}\in \partial (rD))\\
&= 1-\E_x [u(\tau^X_{rD}); X_{\tau^X_{rD}}\in \partial (rD)]\\
&= 1-\E_x [u(\tau^X_{rD})]\\
&\leq 1-\E_x [u(\tau^X_{rD}); \tau^X_{rD}\leq {\rm diam}(rD)^2].
\end{aligned}\end{equation}
Since $u(t)$ is positive and decreasing with $ u(0+)=1$, then for $r\in (0, 1/{\rm diam}(D))$, we have
\begin{equation}\label{e:2.48}\begin{aligned}
\E_x [u(\tau^X_{rD}); \tau^X_{rD}<r^2]&\geq u(1)\P_x [\tau^X_{rD}\leq {\rm diam}(rD)^2]\\
&\geq u(1)\P_x(X_{{\rm diam}(rD)^2}\in (rD)^c)\\
&\geq u(1)\P_x(X_{{\rm diam}(rD)^2}\in B(x, {\rm diam}(rD))^c)\\
&\geq c_3\P_x(W_{{\rm diam}(rD)^2}\in B(x, {\rm diam}(rD))^c)\\
&\geq c_4,
\end{aligned}\end{equation}
where  the fourth inequality is due to \eqref{e:1.3'} and $c_k=c_k(d, \lambda_0)\in (0, 1), k=3, 4.$
Thus, by \eqref{e:2.47} and \eqref{e:2.48},
\begin{equation}\label{e:2.32}
\int_{rD} U^{Z^{rD}}\ast  q_{rD} (x, z)\,dz\leq 1-c_4.
\end{equation}
By Lemma \ref{L:2.5'}, there exists $c_5=c_5(d, \lambda_0)>1$ such that $\sup_{z\in rD}\E_z \tau_{rD}\leq \sup_{z\in rD}\E_z \tau_{B(z, r{\rm diam}(D))} \leq c_5{\rm diam}(rD)^2\leq c_5$ for any $r\in (0, 1/{\rm diam}(D)).$
Hence,
$$\int_{rD}\int_{rD} (U^{Z^{rD}}\ast  q_{rD})^{N+1} (x, z) \E_z \tau_{rD}\,dy\,dz\leq c_5(1-c_4)^{N+1}\rightarrow 0, \quad N\rightarrow\infty.$$
Hence, by this together with \eqref{e:2.21},
\begin{equation}\label{e:2.22}
\E_x \tau_{rD}=\E_x \zeta+\sum_{n=1}^\infty \int_{rD} (U^{Z^{rD}}\ast  q_{rD})^n (x, z) \E_z \zeta\,dz, \quad x\in rD.
\end{equation}

It follows from  Proposition \ref{P:2.1} that there exists $c_6=c_6(d, \lambda_0, \ell, \phi, R_0, \Lambda_0, {\rm diam}(D))$ such that $U^{Z^{rD}}(x, y)\leq c_6G^W_{rD}(x, y)$ for $x, y\in rD$ and $r\in (0, 1).$
By \eqref{e:2.6} and a simple calculation, there exists $c_7=c_7(d, R_0, \Lambda_0)>1$ such that for $x\in rD,$
\begin{equation}\label{e:2.45}
c^{-1}_7\delta_{rD}(x){\rm diam}(rD)\leq \E_x \tau^W_{rD}\leq c_7\delta_{rD}(x){\rm diam}(rD).
\end{equation}
Hence, for any $x\in rD,$
 $$\E_x \zeta=\int_{rD}U^{Z^{rD}}(x, y)\,dy\leq c_6\int_{rD}G^W_{rD}(x, y)\,dy = c_6\E_x \tau^W_{rD}\leq c_6c_7\delta_{rD}(x) {\rm diam}(rD).$$
Thus we have by \eqref{e:2.22},
$$\begin{aligned}
\E_x\tau_{rD}\leq c_6c_7\delta_{rD}(x){\rm diam}(rD) +c_6c_7{\rm diam}(rD) \sum_{n=1}^\infty\int_{rD} (c_6G^W_{rD} \ast  q_{rD})^n (x, z) \delta_{rD}(z) \,dz.
\end{aligned}$$
By applying Lemma \ref{L:2.4} and the induction, there exists $c_8=c_8(d, \lambda_0, \ell, \phi, R_0, \Lambda_0, {\rm diam}(D))$ such that
$$\int_{rD} (G^W_{rD}\ast  q_{rD})^n (x, z)\delta_{rD}(z) \,dz\leq (c_8\Psi({\rm diam}(rD)))^n\delta_{rD}(x),$$
where $\Psi(r)\rightarrow 0$ as $r\rightarrow 0.$
Hence,
$$
\E_x \tau_{rD}\leq c_6c_7\delta_{rD}(x){\rm diam}(rD) +c_6c_7\delta_{rD}(x){\rm diam}(rD) \sum_{n=1}^\infty (c_6c_8\Psi(r{\rm diam}(D)))^n.
$$
Let $\delta_1\in (0, R_0)$ be a small constant such that
 $c_9:=\sum_{n=1}^\infty (c_6c_8\sup_{s\leq \delta_1}\Psi(s))^n<\infty.$ Hence, for $r<\delta_1/{\rm diam}(D),$ we have
$$\E_x \tau_{rD}\leq c_6c_7(1+c_9)\delta_{rD}(x){\rm diam}(rD).$$
Thus it follows from \eqref{e:2.45} that $$\E_x \tau_{rD}\leq c_6c_7^2(1+c_9)\E_x \tau^W_{rD}.$$
This together with \eqref{e:2.19'} yields the desired conclusion.
\qed

\section{Exit distribution estimates}

In this section, we shall establish the exit distribution estimates for $Y$ from a $C^{1, \alpha}$ open set in Proposition \ref{P:3.9}.
When $d\geq 2,$ we mainly use the "box" method developed by Bass and Burdzy in \cite{Ba, BB1}.

We say an open set $D\subset \R^d$ is Greenian with respect to $Y$ if the Green function $G_{D}(x, y)$ of $Y$ in $D$ exists and is not identically infinite.
For any  Greenian  (with respect to $Y$)  open set $D$ in $\R^d$, and for any Borel subset $A\subset D,$ we define
\begin{eqnarray}
 {\rm Cap}_D(A) &:=& \sup\big\{\mu(A): \mu \: \mbox{is a measure supported on} \: A    \nonumber \\
 && \hskip 0.8truein
 \mbox{with} \: \sup_{x\in D}\int_D G_D(x, y)\mu(dy)\leq 1\big\}. \label{e:3.0}
\end{eqnarray}
 The following facts are known; see \cite{CF, FOT}.  Every function $u\in W^{1, 2}(\R^d)$ has an $ \mathcal{E}$-quasi-continuous version,
which is unique $ \mathcal{E}$-quasi-everywhere ($ \mathcal{E}$-q.e. in abbreviation) on $\R^d$.
We always represent $u\in W^{1, 2}(\R^d)$ by its $ \mathcal{E}$-quasi-continuous version.
For a  Greenian open set $D$ and $A\subset D$,
 \begin{equation}\label{e:3.1}
{\rm Cap}_D(A)
=\inf\big\{\mathcal{E}(u, u): u\in W^{1, 2}(\R^d),  \  u\geq 1 \  \mathcal{E} \hbox{-q.e. on }   A
\hbox{ and }  u=0 \     \mathcal{E} \hbox{-q.e. on }   D^c   \big\}.
\end{equation}
We use ${\rm Cap}^X_D(\cdot)$ and ${\rm Cap}^W_D(\cdot)$ to denote the capacity measure of diffusion process $X$ and Brownian motion $W$ in $D.$
Recall that $\mathcal{E}_X$ is the Dirichlet form of $X.$ Let $\mathcal{E}_W$ be the Dirichlet form of $W.$
It follows from the uniform ellipticity \eqref{e:1.2} of $X,$ $\mathcal{E}_X\geq\lambda^{-1}_0 \mathcal{E}_W.$
 Since $\mathcal{E}_X\leq \mathcal{E}$ by \eqref{e:1.6'},  for any Greenian open set $D\subset \R^d,$
  \begin{equation}\label{e:3.2}
\lambda^{-1}_0 {\rm Cap}^W_D(A)\leq  {\rm Cap}^X_D(A)\leq {\rm Cap}_D(A) \quad \mbox{for every Borel subset} \: A\subset D.
  \end{equation}

\begin{defn} \rm
Suppose $U$ is an open set in $\R^d.$ A real-valued function $u$ defined on $\R^d$ is said to be
 harmonic in   $U$ with respect to $Y$
if for every open set $B$  whose closure is a compact subset of $U,$
   $$\E_x |u(Y_{\tau_B})|<\infty \quad \mbox{and} \quad
 u(x)=\E_x u(Y_{\tau_B}) \quad \mbox{for  each} \quad  x\in B.$$
  In particular, we say $u$  is   regular harmonic in   $U$ with respect to $Y$
if    $\E_x |u(Y_{\tau_U})|<\infty$ and
$u(x)=\E_x u(Y_{\tau_U})$  for each  $x\in U.$
\end{defn}

\begin{prp}[Harnack inequality]\label{P:3.3}
 Let $d\geq 2.$
Suppose $H$ satisfies $L^a(\gamma, c_L)$ and $U^a(\delta, C_U)$ with  $\delta<2$ for some $a>0.$
 There exists $C=C(d, \lambda_0, \ell,  \phi)$ such that for any $x_0\in\R^d, r\in (0, 1)$ and nonnegative harmonic function $h$ in $B(x_0, 2r)$  with respect to $Y,$
$$h(x)\leq C h(y), \quad x, y\in B(x_0, r).$$
\end{prp}

\pf We first consider the case $d\geq 3.$
By Proposition \ref{P:2.1}, there exists $c_1=c_1(d, \lambda_0, \ell, \phi)$ such that for any $r\in (0, 1)$ and $x, y\in B(x_0, r)$
$$G_{B(x_0, 2r)}(x, y)\geq U^{Z^{B(x_0, 2r)}(x, y)}\geq c_1 r^{2-d}=c_2\dfrac{1}{{\rm Cap}^W_{B(x_0, 2r)}(B(x_0, r))}.$$
Then by \eqref{e:3.2}, a similar argument in \cite[Lemma 4.1]{RSV}, there exists $c_3=c_3(d, \lambda_0, \ell, \phi)$ such that for any $x_0\in\R^d, r\in (0, a^{-1/2}/4)$ and any closed subset $A$ of $B(x_0, r),$
\begin{equation}\label{e:2.33n}
\P_y(T_A< \tau_{B(x_0, 2r)})\geq c_3\dfrac{{\rm Cap}^W_{B(x_0, 2r)}(A)}{{\rm Cap}^W_{B(x_0, 2r)}(B(x_0, r)) }, \quad y\in B(x_0, r).
\end{equation}
Suppose $H$ satisfies $L^a(\gamma, c_L)$ and $U^a(\delta, C_U)$ with  $\delta<2$ for some $a>0.$
Then $\psi$ is a non-negative function  satisfying $L_{a^{-1/2}}(2\gamma, c_L)$ and $U_{a^{-1/2}}(2\delta, C_U).$
By \eqref{e:2.9},
there exists $c_4=c_4(d, \phi)$ such that
\begin{equation}\label{e:3.8}
j(u)\leq c_4j(2u) \quad \mbox{for} \quad u\in (0, a^{-1/2}/2).
\end{equation}
Since $S_t$ is a complete subordinator, it follows from \cite[Lemma 2.1]{KSV3} that there is a positive constant
$c_5$ such that $\mu(t)\leq c_5 \mu(t+1)$ for $t\geq a^{-1/2}.$ Then by \cite[(2.7)]{KSV2}, there exists
 $c_6=c_6(d, \phi)$ such that
\begin{equation}\label{e:3.9}
 \quad j(u)\leq c_6j(u+1)\quad \mbox{for} \quad u\geq a^{-1/2}.
\end{equation}
Hence, by \eqref{e:2.33n}-\eqref{e:3.9}, Lemma \ref{L:2.5'} and a very similar argument in \cite[Theorem 4.5]{RSV},
the desired conclusion is obtained when $d\geq 3.$.

When $d\geq 2,$ the conclusion is obtained by a similar argument in \cite[Proposition 2.2]{KSV2}.
\qed

\medskip

\begin{lem}\label{L:3.2}
Suppose $H$ satisfies $L^a(\gamma, c_L)$ and $U^a(\delta, C_U)$ with  $\delta<2$ for some $a>0.$
There exists $C=C(d, \lambda_0, \phi)>0$ such that for any $r\in (0, 1)$ and $x\in B(0,  r)\setminus \{0\},$
$$
G_{B(0, r)}(x, 0 )\leq \left\{\begin{array}{ll}
C|x|^{2-d}, & \quad d\geq 3\\
C\log(3r/|x|),  & \quad d=2.
\end{array}\right.$$
\end{lem}

\pf The proof is similar to \cite[Lemma 4.6]{CKSV1}.
Fix $x\in  B(0, r)\setminus \{0\}$ and let $\rho:=|x|/(3r).$ Then $\rho\in (0, 1/4).$ Since $\overline{B(0, \rho r)}= \overline{B(0, |x|/3)}$ is a compact subset of $B(0, r),$ there exists a capacitary measure $\mu_a$ for $\overline{B(0, \rho r)}$ such that ${\rm Cap}_{B(0, r)}(\overline{B(0, \rho r)})=\mu_a(\overline{B(0, \rho r)}).$
Note that $y\mapsto G_{B(0, r)}(x, y)$ is harmonic with respect to $Y$ in $B(0, 2\rho r)=B(0, 2|x|/3).$ By the uniform Harnack inequality in Proposition \ref{P:3.3}, we have
$$\begin{aligned}
1&\geq \int_{\overline{B(0, \rho r)}} G_{B(0, r)}(x, y)\mu_\rho(dy)\\
&\geq \left(\inf_{y\in \overline{B(0, \rho r)}} G_{B(0, r)}(x, y)\right)\mu_\rho(\overline{B(0, \rho r)})\\
&\geq c_1G_{B(0, r)}(x, 0) {\rm Cap}_{B(0, r)}(\overline{B(0, \rho r)})\\
&\geq c_1\lambda^{-1}_0 G_{B(0, r)}(x, 0) {\rm Cap}^W_{B(0, r)}(\overline{B(0, \rho r)})\\
\end{aligned}$$
where the constant $c_1$ is independent of $r\in (0, a^{-1/2}/4)$ and we used \eqref{e:3.2} in the last inequality.
Hence,
\begin{equation}\label{e:3.6'}
G_{B(0, r)}(x, 0)\leq  \dfrac{c_1^{-1}\lambda_0}{{\rm Cap}^W_{B(0, r)}(\overline{B(0, \rho r)})}=\dfrac{c_1^{-1}\lambda_0}{{\rm Cap}^W_{B(0, r)}(\overline{B(0, |x|/3)})}.
\end{equation}
It is known that (see e.g. \cite[Lemma 4.5]{CKSV1}) there exists $c_2>0$ such that for any $\rho\in (0, 1/4),$
\begin{equation}
{\rm Cap}^W_{B(0, 1)}(\overline{B(0, \rho)})\geq
\left\{\begin{array}{ll}
c_2|\rho|^{d-2}, & \quad d\geq 3\\
c_2/\log(1/|\rho|),  & \quad d=2.
\end{array}\right.
\end{equation}
By the scaling property of $W$,
 we have $ G^W_{B(0, 1)}(x, y)=r^{d-2}G^W_{B(0, r)}(rx, ry)$ for $x, y\in B(0, 1).$ Hence, it follows from \eqref{e:3.0} that
\begin{equation}\label{e:3.8'}
r^{d-2}{\rm Cap}^W_{B(0, 1)}(\overline{B(0, \rho)})
 ={\rm Cap}^W_{B(0, r)}(\overline{B(0, \rho r)}).
\end{equation}
The conclusion now follows from \eqref{e:3.6'}-\eqref{e:3.8'}.
\qed

\medskip

Recall that $\delta_1$ is the constant in Proposition \ref{L:2.5}.

\begin{lem}\label{L:3.4}
Let $d\geq 2.$
Suppose $H$ satisfies the assumption (A1).
There exists $C=C(d, \lambda_0, \ell, \phi)$ such that for any  $x_0\in\R^d$ and $r\in (0, \delta_1),$
$$G_{B(x_0, r)}(x, y)\leq CG^W_{B(x_0, r)}(x, y) \quad {\rm for} \quad x\in B(x_0, r/4) \quad {\rm and}\quad y\in B(x_0, r)\setminus B(x_0, r/2).$$
\end{lem}

\pf Let $x_0\in\R^d$ and $r\in (0, \delta_1).$ For the simplicity of notation, let $B_r:=B(x_0, r).$
Since for each $y\in B_r\setminus B_{r/2},$ $x\mapsto G_{B_r}(x, y)$ is harmonic with respect to $Y$ in $B_{r/2},$ by the Harnack inequality Proposition \ref{P:3.3}, there exists $c_1=c_1(d, \lambda_0, \phi)$ such that for each $x\in B_{r/4},$
$$G_{B_r}(x, y)\leq \dfrac{c_1}{r^d}\int_{B_{r/4}} G_{B_r}(u, y)\,du\leq c_1r^{-d} \E_y \tau_{B_r}, \quad y\in B_r\setminus B_{r/2}.$$
By Proposition \ref{L:2.5} and \eqref{e:2.45}, there exist $c_k=c_k(d, \lambda_0, \ell, \phi), k=2, 3$ such that for  $r\in (0, \delta_1),$
$$\E_y \tau_{B_r}\leq c_2\E_y \tau^W_{B_r} \leq c_3r\delta_{B_r}(y), \quad  y\in B_r\setminus B_{r/2}.$$
Hence, $$G_{B_r}(x, y)\leq c_1c_3r^{1-d}\delta_{B_r}(y) \quad {\rm for} \quad x\in B_{r/4} \quad {\rm and}\quad y\in B_r\setminus B_{r/2}.$$
By  \eqref{e:2.6}, when  $d\geq 3,$
$$G_{B_r}(x, y)\leq c_1c_3r^{1-d}\delta_{B_r}(y)\leq c_4G^W_{B_r}(x, y) \quad {\rm for} \quad x\in B_{r/4} \quad {\rm and}\quad y\in B_r\setminus B_{r/2}.$$
When $d=2,$ note that $\log(1+s)\geq c_5 s$ for $s\in (0, 1),$ then by  \eqref{e:2.6},  we have
$$G_{B_r}(x, y)\leq c_1c_3r^{-1}\delta_{B_r}(y) \leq c_5 \log(1+\delta_{B_r}(y)/r)\leq c_6G^W_{B_r}(x, y) \quad {\rm for} \quad x\in B_{r/4} \quad {\rm and}\quad y\in B_r\setminus B_{r/2}.$$
\qed

\medskip

Let $D$ be a $C^{1, \alpha}$ open set  with characteristics $(R_0, \Lambda_0).$
For each $\lambda\geq 1,$ since $\Gamma_\lambda(x):=\lambda \Gamma(x/\lambda)$ is the graph function of the boundary of $\lambda D,$
it is easy to see that $\lambda D$ is  a $C^{1, \alpha}$ open set  with characteristics $(\lambda R_0, \Lambda_0).$
By a similar argument in \cite[Lemma 2.2]{S} for $C^{1, 1}$ open sets, for each $r\in (0, R_0/2),$ there exists $L=L(R_0, \Lambda_0, d)>1$
such that for any $z\in\partial D,$ there is a $C^{1,\alpha}$ connected open set  $U_{z, r}\subset D$  such that $D\cap B(z, r)\subset U_{z, r}\subset D\cap B(z, 2r)$ and  $r^{-1}U_{z, r}$  is a $C^{1, \alpha}$ open set with characteristics $(R_0/L, L\Lambda_0).$
Hence, for each $r\in (0, R_0/2),$ $U_{z, r}$ is a $C^{1, \alpha}$ open set  with characteristics $(rR_0/L, L\Lambda_0/r^\alpha).$
In the following, we always use $U_{z, r}$ to denote such $C^{1,\alpha}$ open set.

Recall that an open set $D$ in $\R^d$ (when $d\geq 2$) is said to be Lipschitz if there exist a localization radius
$R_0 >0$ and a constant $\Lambda_0>0$ such that for every $z\in\partial D,$ there exist a Lipschitz function
$\Gamma=\Gamma_z: \R^{d-1}\rightarrow \R$ satisfying
$$
\Gamma(0)=\nabla \Gamma(0)=0,
|\Gamma(x)-\Gamma(y)|\leq \Lambda_0 |x-y|,
$$
 and an orthonormal coordinate system
$CS_z: y=(y_1, \cdots, y_{d-1}, y_d)=:(\wt {y}, y_d)\in \R^{d-1}\times \R$ with its origin at $z$ such that
$$
B(z, R_0)\cap D=\{y=(\wt {y}, y_d)\in B(0,R_0) \hbox{ in } CS_z: y_d>\Gamma(\wt {y})\}.
$$
The pair $(R_0, \Lambda_0)$ is called the characteristics of the Lipschitz open set $D.$

Suppose $D$ is a $C^{1, \alpha}$ open set  with characteristics $(R_0, \Lambda_0).$
Then $D$ is a Lipschitz open set with characteristics $(R_0, \Lambda_0).$
 It is well known that
there exists   $ \kappa=\kappa(R_0, \Lambda_0)\in (0, 1/4)$ such that for $r\in (0, R_0)$
and $z\in\partial D,$
\begin{equation}\label{e:3.10}
\hbox{there exists } z_r\in D\cap \partial B(z, r)  \hbox{ with }
\kappa r\leq \delta_D(z_r)< r.
 \end{equation}
   In the following,  we always use $\kappa $ to denote the positive constant   in \eqref{e:3.10}.

 \begin{lem}\label{L:2.6'}
Let $D$ be a $C^{1, \alpha}$ open set with characteristics $(R_0, \Lambda_0).$
There exists a positive constant $C=C(d, \lambda_0, \ell, \phi, R_0, \Lambda_0)\in (0, 1)$ such that for any $z_0\in\partial D$ and $r\in (0, R_0/4),$
$$
\P_x(Y_{\tau_{D_r(z_0)}}\in  D\cap \partial B(z_0, r))\geq C\delta_{D}(x)/r, \quad x\in D_{\kappa r/2}(z_0),
$$
where $D_r(z_0):=D\cap B(z_0, r).$
\end{lem}

\pf Recall that the potential density function $u$ of $S_t$ is  strictly positive and decreasing continuous   on $[0, \infty)$ with $u(0+)=1$.

Let $z_0\in \partial D$ and $r\in (0, R_0/4).$ Let $Z^{D_r(z_0)}_t:=X^{D_r(z_0)}(S_t)$ be the subordinate killed diffusion in $D_r(z_0).$
 We will use $\zeta$ to denote the life time of the process $Z^{D_r(z_0)}_t.$
 Let $x\in D_{\kappa r/2}(z_0).$
 By  Theorem \ref{T:2.2}, we have
\begin{equation}\label{e:3.10}\begin{aligned}
&\P_x(Y_{\tau_{D_r(z_0))}}\in  D\cap \partial B(z_0, r))\\
&\geq  \P_x(Z^{D_r(z_0)}_{\zeta-}\in  D\cap \partial  B(z_0, r))\\
&=   \E_x [u(\tau^X_{D_r(z_0)}); X_{\tau^X_{D_r(z_0)}}\in D\cap \partial  B(z_0, r)]\\
&\geq   \E_x [u(\tau^X_{D_r(z_0)});  \tau^X_{D_r(z_0)}\leq r^2t \  \hbox{ and }  \  X_{\tau^X_{D_r(z_0)}}\in D\cap \partial  B(z_0, r)]\\
&\geq   \inf_{s\in (0, t)} u(s)\cdot \P_x(\tau^X_{D_r(z_0)}\leq r^2t \  \hbox{ and }  \  X_{\tau^X_{D_r(z_0)}}\in D\cap \partial  B(z_0, r))\\
&\geq   u(t)\cdot  \left( \P_x(X_{\tau^X_{D_r(z_0)}}\in D\cap \partial  B(z_0, r))-\P_x(\tau^X_{D_r(z_0)}>r^2 t)  \right) .
\end{aligned}\end{equation}
Let $U_{z_0, 2r}$ be a  $C^{1, \alpha}$ domain with characteristics
$(2rR_0/L, \Lambda_0 L/(2r)^\alpha)$ such that $D_{2r}(z_0)\subset U_{z_0, 2r}\subset D_{4r}(z_0).$
Let $x_0$ be a point in $D\cap \partial B(z_0, r/2)$ such that $\delta_{D}(x_0)\geq \kappa r/2.$  Let  $y_0$ be a point in $D\cap \partial B(z_0, 3r/2)$  such that  $\delta_{D}(y_0)\geq 3\kappa r/2.$
Note that $G^X_{U_{z_0, 2r}}(\cdot, y_0)$ is harmonic in $D_r(z_0).$
By the scale invariant boundary Harnack principle for $X$ on Lipschitz domain (see e.g. \cite{CFMS}) and \eqref{e:2.6},
there exist $c_k=c_k(d, \lambda_0, \ell, R_0, \Lambda_0)\in (0, 1), k=1, 2$ such that for $x\in D_{\kappa r/2}(z_0),$
\begin{equation}\label{e:3.12'}
\P_x(X_{\tau^X_{D_r(z_0)}}\in D\cap \partial  B(z_0, r))\geq c_1\dfrac{G^X_{U_{z_0, 2r}}(x, y_0)}{G^X_{U_{z_0, 2r}}(x_0, y_0)}\geq c_2\dfrac{\delta_{D}(x)}{r}.
\end{equation}
On the other hand, by Proposition \ref{L:2.5} and \eqref{e:2.45},  there exist $c_k=c_k(d, \lambda_0, \ell, \phi, R_0, \Lambda_0)>0, k=3, 4$ such that  for $x\in D_{\kappa r/2}(z_0),$
\begin{equation}\label{e:3.13'}
\P_x \big(\tau^X_{D_r(z_0)}>r^2 t \big)\leq \dfrac{\E_x\tau^X_{U_{z_0, 2r}}}{r^2t}\leq c_3\dfrac{\E_x\tau^W_{U_{z_0, 2r}}}{r^2t}\leq c_4\dfrac{\delta_{D}(x)}{rt}.
\end{equation}
We choose a large enough constant  $t_0=t_0(R_0, \Lambda_0)>0$ such that $c_2-c_4/t_0 \geq \frac{c_2}{2}.$
Hence, by \eqref{e:3.10}-\eqref{e:3.13'},
$$\P_x(Y_{\tau_{D_r(z_0))}}\in  D\cap \partial B(z_0, r))\geq u(t_0)\dfrac{c_2}{2}\dfrac{\delta_{D}(x)}{r}, \quad x\in D_{\kappa r/2}(z_0).$$
  This proves the lemma.
\qed

\medskip

The following Lemma follows from a similar argument as \cite[Lemma 4.1]{CKSV} or \cite[Lemma 5.1]{KSV2} with the process $X$ in place of $W$, we omit the proof here.

\begin{lem}\label{L:2.6}
Let $D$ be a Lipschitz open set with characteristics $(R_0, \Lambda_0).$
There exists a positive constant $\rho_0=\rho_0(d, \lambda_0, \ell, \phi, R_0, \Lambda_0)\in (0, 1)$ such that for any $ x\in D$ with $\delta_{D}(x)\leq R_0/2,$
$$\P_x(Y_{\tau_{B(x, 2\delta_{D}(x))\cap D}}\in D^c)\geq \rho_0.
$$
\end{lem}

\medskip
Let $D$ be a $C^{1, \alpha}$ open set in $\R^d$ with characteristics $(R_0, \Lambda_0).$
Let $x\in  D.$ Let $z_x\in \partial D$ such that $|x-z_x|=\delta_D(x).$
Let $y=(\tilde y, y_{d})\in \R^{d-1}\times \R$ be the coordinate in $CS_{z_x},$
 define $$\rho_\Gamma(y):=y_{d}-\Gamma(\tilde y).$$
  For $r\in (0, R_0),$ we define the "box":
$$\Delta(x, a, r):=\{y\in D \: \mbox{in} \: CS_{z_x}: 0<\rho_\Gamma(y)<a\}\cap B(z_x, r).$$

\begin{lem}\label{L:3.5}
 Suppose $d\geq 2$ and $H$ satisfies the assumption (A1).
Suppose $D$ is a $C^{1, \alpha}$ open set with characteristics $(R_0, \Lambda_0)$ in $\R^d$. For each $M> 1,$ there exists a constant   $C_M=C_M(d, \lambda_0, \ell, \phi, R_0, M)$  such that for any
$ s\in (0, (R_0\wedge a^{-1/2})/(2M)]$ and any $x\in D$ with $\rho_\Gamma(x)<s$ in $CS_{z_x},$
$$\P_x(Y_{\tau_{\Delta(x, s, Ms)}}\in \Delta(x, s, 2Ms)\setminus \Delta(x, s, Ms))\leq C_M,$$
where $C_M\rightarrow 0$ as $M\rightarrow\infty.$
\end{lem}

\pf Let $M\geq 1$ and $ s\in (0, (R_0\wedge a^{-1/2})/(2M)].$
Without loss of generality, we assume  the constant  $a=1$ in the assumption (A1).
Then $R_0\wedge a^{-1/2}=R_0.$
 Let $x\in D$ with $\rho_\Gamma(x)<s$ in $CS_{z_x}.$
For the simplicity of notation, for each $r>0,$ let $B_r:=B(x, r).$ Let $y$ be a fixed point on $\Delta(x, s, (M+1)s)\setminus \Delta(x, s, Ms).$
Note that $G_{B_{3Ms}}(\cdot, y)$  is regular harmonic in $B_{Ms}.$ Hence,
\begin{equation}\label{e:3.11'}\begin{aligned}
G_{B_{3Ms}}(x, y)&=\E_x G_{B_{3Ms}}(Y_{\tau_{B_{Ms}}}, y)\\
&\geq \E_x [G_{B_{3Ms}}(Y_{\tau_{B_{Ms}}}, y); Y_{\tau_{B_{Ms}}}\in \Delta(x, s, (M+1)s)\setminus \Delta(x, s, Ms)]\\
&\geq \inf_{z\in \Delta(x, s, (M+1)s)\setminus \Delta(x, s, Ms)}G_{B_{3Ms}}(z, y)\P_x( Y_{\tau_{B_{Ms}}}\in \Delta(x, s, (M+1)s)\setminus \Delta(x, s, Ms)).
\end{aligned}\end{equation}
Let $Z^{B_{3Ms}}_t:=X^{B_{3Ms}}(S_t)$ and  $U^{Z^{B_{3Ms}}}$ be the Green function of the subordinate killed process $Z^{B_{3Ms}}.$
By Proposition \ref{P:2.1} and \eqref{e:2.5'}, there exist $c_k=c_k(d, \lambda_0, \ell, \phi)>0, k=1, 2$ such that for any $z\in \Delta(x, s, (M+1)s)\setminus \Delta(x, s, Ms),$
$$
G_{B_{3Ms}}(z, y)\geq U^{Z^{B_{3Ms}}}(z, y)\geq c_1G^W_{B_{3Ms}}(z, y)\geq c_2g_{B_{3Ms}}(z, y).
$$
Hence, there exists $c_3=c_3(d, \lambda_0,  \phi)>0$ such that for any $z\in \Delta(x, s, (M+1)s)\setminus \Delta(x, s, Ms),$
\begin{equation}\label{e:3.11}
G_{B_{3Ms}}(z, y)\geq c_2g_{B_{3Ms}}(z, y)\geq
\left\{\begin{array}{ll}
c_3s^{2-d}, & \quad d\geq 3\\
c_3\log(M),  & \quad d=2.
\end{array}\right.
\end{equation}
By Lemma \ref{L:3.2}, there exists $c_4=c_4(d, \lambda_0,  \phi)>0$ such that
\begin{equation}\label{e:3.12}
G_{B_{3Ms}}(x, y)\leq
\left\{\begin{array}{ll}
c_4(Ms)^{2-d}, & \quad d\geq 3\\
c_4\log 3,  & \quad d=2.
\end{array}\right.
\end{equation}
Hence, when $d\geq 2,$ by \eqref{e:3.11'}-\eqref{e:3.12},
\begin{equation}\label{e:3.7}
\P_x(Y_{\tau_{\Delta(x, s, Ms)}}\in \Delta(x, s, (M+1)s)\setminus \Delta(x, s, Ms))\leq c_5(\log M)^{-1}.
\end{equation}

For the simplicity of notation,  let $\Delta_2:=\Delta(x, s, 2Ms).$
By the L\'evy system formula of $Y$,
\begin{equation}\label{e:2.25}\begin{aligned}
&\P_x(Y_{\tau_{\Delta(x, s, Ms)}}\in \Delta(x, s, 2Ms)\setminus \Delta(x, s, (M+1)s))\\
&=\P_x(Y_{\tau_{\Delta(x, s,  Ms)}}\in  \Delta_2 \cap \overline B^c(x, (M+1)s))\\
&=\int_{\overline B^c(x, (M+1)s)\cap  \Delta_2}\int_{\Delta(x,s,  Ms)}G_{\Delta(x,s,  Ms)}(x, y)J(y, z)\,dydz\\
&=\int_{\overline B^c(x, (M+1)s)\cap  \Delta_2}\int_{\Delta(x,s,  Ms/2)}G_{\Delta(x,s,  Ms)}(x, y)J(y, z)\,dydz\\
&+\int_{\overline B^c(x, (M+1)s)\cap  \Delta_2}\int_{\Delta(x,s,  Ms)\setminus \Delta(x,s,  Ms/2)}G_{\Delta(x,s,  Ms)}(x, y)J(y, z)\,dydz\\
&=:I+II.
\end{aligned}\end{equation}
For the first term,  by \eqref{e:2.9'} and \eqref{e:2.4}, $J(y, z)\leq c_6\frac{\psi_0(|y-z|)}{|y-z|^{d+2}}\leq c_7|y-z|^{-(d+2)}$ for $y, z\in\R^d$ with $|y-z|\leq 1.$ Then by Lemma \ref{L:2.5'} and  $Ms\leq 1/2$, we have
\begin{equation}\label{e:3.13}\begin{aligned}
I&\leq c_7\E_x \tau_{\Delta(x, s,  Ms)}\sup_{y\in \Delta(x,s,  Ms/2)}\int_{\Delta_2\setminus B(x, (M+1)s)} |y-z|^{-(d+2)}\,dz\\
&\leq c_7\E_x \tau_{B(x, Ms)}\int_{|\tilde z|\geq Ms/2}
\int_{\Gamma(\tilde z)}^{\Gamma(\tilde z)+s}|\tilde z|^{-(d+2)} \,dz_{d}\, d\tilde z\\
&\leq c_8(Ms)^2 (Ms)^{-3}s\\
&\leq c_8M^{-1}.
\end{aligned}\end{equation}

Let $W$ be a Brownian motion independent of $S_t.$ Let $Y^0_t:=W_{S_t}.$
 For each open set $B,$ denote by $G^0_B(x, y)$  the Green function of $Y^0$ in $B$ and $\tau^0_B$ the first exit time for $Y^0$ from $B,$ respectively.
It follows from \cite[Corollary 1.7]{BK} that $G^0_{B(0, 1)}(x, y)\asymp G^W_{B(0, 1)}(x, y)$ for $x, y\in B(0, 1).$
For the second term in \eqref{e:2.25}, let $\bar x=(Ms)^{-1}x.$  Note that by \eqref{e:2.9'}, $J(y, z)\asymp \frac{\psi_0(|y-z|)}{|y-z|^{d+2}}$ for $|y-z|\leq 4.$ By Lemma \ref{L:3.4}, we have
\begin{equation}\label{e:3.14}\begin{aligned}
II&\leq c_9\int_{\overline B^c(x, (M+1)s)\cap \Delta_2}\int_{B(x, Ms)\setminus B(x, Ms/2)}G^W_{B(x, Ms)}(x, y)J(y, z)\,dydz\\
&\leq c_{10}\int_{\overline B^c(x, (M+1)s)\cap \Delta_2}\int_{B(x, Ms)\setminus B(x, Ms/2)}G^W_{B(x, Ms)}(x, y)\dfrac{\psi_0(|y-z|)}{|y-z|^{d+2}}\,dydz\\
&= c_{10}\int_{\overline B^c(x, (M+1)s)\cap \Delta_2}\int_{B(x, Ms)\setminus B(x, Ms/2)}(Ms)^{2-d}G^W_{B(\bar x, 1)}(\bar x, (Ms)^{-1}y)\dfrac{\psi_0(|y-z|)}{|y-z|^{d+2}}\,dydz\\
&= c_{10}\int_{\overline B^c(\bar x, \frac{M+1}{M})\cap (Ms)^{-1}\Delta_2}\int_{B(\bar x, 1)\setminus B(\bar x, 1/2)}G^W_{B(\bar x, 1)}
(\bar x, y)\dfrac{\psi_0(|Ms(y-z)|)}{|y-z|^{d+2}}\,dydz\\
&\leq  c_{11}\int_{\overline B^c(\bar x, 1)\cap (Ms)^{-1}\Delta_2}\int_{B(\bar x, 1)}
G^W_{B(\bar x, 1)}(\bar x,y)\dfrac{\psi_0(|y-z|)}{|y-z|^{d+2}}\,dydz\\
&\leq  c_{12}\int_{\overline B^c(\bar x, 1)\cap (Ms)^{-1}\Delta_2}\int_{B(\bar x, 1)}
G^0_{B(\bar x, 1)}(\bar x,y)j(|y-z|)\,dydz\\
&= c_{12}\P_{\bar x}\left(Y^0_{\tau^0_{B(\bar x, 1)}}\in \overline B^c(\bar x, 1)\cap (Ms)^{-1}\Delta_2; \: Y^0_{\tau^0_{B(\bar x, 1)}-}\neq Y^0_{\tau^0_{B(\bar x, 1)}}\right)\\
&=:c_M,
\end{aligned}\end{equation}
where in the third line we used the scaling property of $ G^W_{B(x, r)}(x, y)=r^{2-d}G^W_{B(r^{-1}x, 1)}(r^{-1}x, r^{-1}y)$ for $r>0,$
and in the fifth line  we used \eqref{e:2.9'} and \eqref{e:2.4}.

Note that $(Ms)^{-1}\Delta_2=\Delta(\bar x, M^{-1}, 2)\rightarrow \varnothing$ as $M\rightarrow\infty.$
 Hence,  $c_M=c_M(d, \lambda_0, \ell, \phi)$ holds for any $ s\in (0, (R_0\wedge a^{-1/2})/(2M)]$  and $c_M\rightarrow 0$ as $M\rightarrow\infty.$
Consequently, by \eqref{e:3.7}-\eqref{e:3.14}, the desired conclusion is obtained.

\qed

\begin{lem}\label{L:2.7}
 Suppose $d\geq 2$ and $H$ satisfies the assumption (A1).
Suppose  $D$ is a $C^{1, \alpha}$ open set with characteristics $(R_0, \Lambda_0)$ in $\R^d.$ There exists a constant   $C=C(d, \lambda_0, \ell, \phi, R_0)$  such that for any
$ M\geq 1,  s\in (0, (R_0\wedge a^{-1/2})/(8M)]$ and any $x\in D$ with  $\rho_\Gamma(x)<s$ in $CS_{z_x},$
$$\P_x(Y_{\tau_{\Delta(x, s, Ms)}}\in \Delta(x, s, (a^{-1/2}\wedge R_0)/2))\leq CM^{-3}.$$
\end{lem}

\pf The proof  combines Lemma \ref{L:3.5} and the  method in \cite[Lemma 3]{BB}.
Without loss of generality, we assume  the constant  $a=1$ in the assumption (A1).
Let $M_0\geq 2$ be a positive integer which will be chosen later and $s\in (0, R_0/(8M_0)).$
 Let $K_0$ be the integer part of $R_0/(2M_0s).$
Then $K_0\geq 4.$
For $i=1,2, \cdots$ and $y\in D,$
define $\Delta_i(y):=\Delta (y, s, iM_0s).$
Note that $\Delta_{K_0}(y)\subset \Delta(y, s, R_0/2).$
Define $\beta_0=1$ and
$$\beta_i:=\sup_{y\in D, \rho_\Gamma(y)<s} \P_y(Y_{\tau_{\Delta_i(y)}}\in \Delta_{K_0}(y)), \quad i\geq 1.$$

By the strong Markov property of $Y,$ for each $i=1, \cdots, K_0-1$  and $y\in D$ with $\rho_\Gamma(y)<s,$
$$\begin{aligned}
&\P_y(Y_{\tau_{\Delta_i(y)}}\in  \Delta_{K_0}(y))=\P_y(Y_{\tau_{\Delta_i(y)}}\in  \Delta_{K_0}(y)\setminus \Delta_i(y))\\
&\leq \P_y(Y_{\tau_{\Delta_1(y)}}\in  \Delta_{K_0}(y)\setminus \Delta_{i-1}(y))
+\sum_{k=2}^{i-1}\P_y(Y_{\tau_{\Delta_1(y)}}\in \Delta_k(y)\setminus\Delta_{k-1}(y), Y_{\tau_{\Delta_i(y)}}\in  \Delta_{K_0}(y))\\
&\leq \P_y(Y_{\tau_{\Delta_1(y)}}\in  \Delta_{K_0}(y) \setminus\Delta_{i-1}(y))\\
&\quad\quad+\sum_{k=2}^{i-1} \P_y(Y_{\tau_{\Delta_1(y)}}\in \Delta_k(y)\setminus\Delta_{k-1}(y))\sup_{u\in \Delta_k(y)\setminus\Delta_{k-1}(y)}\P_u(Y_{\tau_{\Delta_{i-k}}(u)}\in  \Delta_{K_0}(u) ).
\end{aligned}$$
 Hence, for $i=1, \cdots, K_0-1,$
\begin{equation}\label{e:2.24}
\P_y(Y_{\tau_{\Delta_i(y)}}\in  \Delta_{K_0}(y))\leq \P_y(Y_{\tau_{\Delta_1(y)}}\in  \Delta_{K_0}(y) \setminus\Delta_{i-1}(y))+\sum_{k=2}^{i-1} \P_y(Y_{\tau_{\Delta_1(y)}}\in \Delta_k(y)\setminus \Delta_{k-1}(y))\beta_{i-k}.
\end{equation}

By Lemma \ref{L:3.5}, there exists a constant $C_{M_0}=C_{M_0}(d, \lambda_0, \ell, \phi)$ such that for $y\in D$ with $\rho_\Gamma(y)<s,$
 \begin{equation}\label{e:3.15}
 \P_y(Y_{\tau_{\Delta_1}(y)}\in  \Delta_2(y)\setminus \Delta_1(y))\leq C_{M_0}
 \end{equation}
 and $C_{M_0}\rightarrow 0$ as $M_0\rightarrow\infty.$
 By \eqref{e:2.8'}, \eqref{e:2.9'} and \eqref{e:2.4},  there exists $c_1>1$ such that $J(z, u)\leq c_1\frac{\psi_0(|z-u|)}{|z-u|^{d+2}}\leq c_2 |z-u|^{-(d+2)}$ for $|z-u|\leq 1.$ Note that $|z-u|\leq 1$ for $z, u\in \Delta_{K_0}(y).$
Thus for each $3\leq k\leq K_0-1$, by the L\'evy system formula for $Y$ and Lemma \ref{L:2.5'},
 \begin{equation}\label{e:2.26}\begin{aligned}
 &\P_y (Y_{\tau_{\Delta_1(y)}}\in \Delta_{K_0}(y)\setminus\Delta_k(y))\\
= &\E_y \int_0^{\tau_{\Delta_1(y)}}\int_{\Delta_{K_0}(y)\setminus\Delta_k(y)} J(Y_s, z)\,dz\, ds\\
\leq & c_2\E_y \tau_{B(y, M_0 s)}\sup_{u\in \Delta_1(y)}\int_{\Delta_{K_0}(y)\setminus\Delta_k(y)}|z-u|^{-(d+2)} \,dz\\
\leq & c_3(M_0 s)^2 \int_{|\tilde z|\geq kM_0s}\int_{\Gamma(\tilde z)}^{\Gamma(\tilde z)+s}|\tilde z|^{-(d+2)} \,dz_{d}\, d\tilde z\\
\leq &c_4M_0^{-1}k^{-3}.
\end{aligned}\end{equation}
Let $\delta_{M_0}:=C_{M_0}\vee c_4M_0^{-1}.$ Then $\delta_{M_0}\rightarrow 0$ as $M_0\rightarrow \infty.$ Note that $\beta_0=1.$
By \eqref{e:2.24}-\eqref{e:2.26}, we have for  $i=1, \cdots, K_0-1,$
$$
\beta_i\leq \delta_{M_0} \sum_{k=2}^i (k-1)^{-3}\beta_{i-k}.
$$
Note that $\beta_i=0$ for $i\geq K_0.$ Hence,
\begin{equation}\label{e:2.27}
\beta_i\leq \delta_{M_0} \sum_{k=2}^i (k-1)^{-3}\beta_{i-k}, \quad i\geq 1.
\end{equation}

In the following, we assert that there exists $c_0=c_0(d,\lambda_0, \ell, \phi)>1$ such that
\begin{equation}\label{e:2.28}
\beta_i\leq c_0(i+1)^{-3}, \quad  i\geq 1.
\end{equation}
Note that $\beta_0=1$ and $\beta_1\leq 1.$  Obviously \eqref{e:2.28} holds for $i=0, 1.$
For $i\geq 1,$ if \eqref{e:2.28} holds for $i,$ by \eqref{e:2.27} we have
$$\begin{aligned}
(i+2)^3\beta_{i+1}&\leq (i+2)^3\delta_{M_0} \sum_{k=2}^{i+1} (k-1)^{-3}c_0(i-k+1)^{-3}\\
&= \delta_{M_0} \sum_{k=1}^{i} k^{-3}c_0(i-k)^{-3} (i+2)^{3}\\
&= \delta_{M_0} \sum_{1\leq k\leq i/2} k^{-3}c_0\left(\dfrac{i+2}{i-k}\right)^3+\delta_{M_0} \sum_{i/2<k<i} \left(\dfrac{i+2}{k}\right)^3c_0(i-k)^{-3}\\
&\leq 6^3c_0\delta_{M_0}\left(\sum_{1\leq k\leq i/2} k^{-3}+\sum_{i/2<k<i}(i-k)^{-3}\right)\\
&\leq 6^3c_0\delta_{M_0}\sum_{k=1}^\infty k^{-3}
\end{aligned}$$
We take $M_0$ large enough such that
$6^3\delta_{M_0}\sum_{k=1}^\infty k^{-3}\leq 1.$
Consequently, \eqref{e:2.28} holds by induction.

For each $M\geq 1$, let $s\in (0, R_0/(8M)].$
If $1\leq M<M_0,$ we have
$$\P_x(Y_{\tau_{\Delta(x, s, Ms)}}\in \Delta(x, s, R_0/2))\leq 1\leq \dfrac{M_0^3}{M^3}.$$
If $M\geq M_0,$
we write $i$ for the integer part of $M/M_0,$ then $i\geq \frac{1}{2}M/M_0.$
Thus by \eqref{e:2.28},
$$\P_x(Y_{\tau_{\Delta(x, s, Ms)}}\in \Delta(x, s, R_0/2))\leq \P_x(Y_{\tau_{\Delta_i(x)}}\in \Delta_{K_0}(x))\leq c_0(i+1)^{-3}\leq 8c_0i^{-3}\leq 8c_0\dfrac{(2M_0)^3}{M^3}.$$
The proof is complete.
\qed

\begin{lem}\label{L:2.15}
Suppose  $H$ satisfies the assumption (A1).
Let $D$ be a $C^{1, \alpha}$ open set with characteristics $(R_0, \Lambda_0)$ in $\R^d$.
There exists $C=C(d, \lambda_0, \ell, \phi, R_0, \Lambda_0)>0$ such that for each $z_0\in \partial D, r\in (0, R_0/4), j\geq 2$ and an open set
$A\subset D\cap B(z_0, R_0/2)$ with  $\Upsilon:={\rm dist}(\Delta(z_0, 2^{-j}r, r), A)>0,$
$$\P_x(Y_{\tau_{\Delta(z_0, 2^{-j}r, r)}}\in A)\leq C2^{-2j}r^2\Upsilon^{-2}, \quad x\in \Delta(z_0, 2^{-j}r, r).$$
\end{lem}

\pf The proof follows the idea in \cite[Lemma 7]{BB}.
Let $z_0\in \partial D$ and $r\in (0, R_0/4).$ For the simplicity of notation, let $\Delta:=\Delta(z_0, 2^{-j}r, r).$  For  each $x\in \Delta,$  let  $B_x:=B(x, 2^{1-j}r)$ and  $C_x:=\Delta\cap B_x.$  Define
$$p_0(x, A):=\P_x(Y_{\tau_{C_x}}\in A), \quad p_{k+1}(x, A):=\E_x[p_k(Y_{\tau_{C_x}}, A); Y_{\tau_{C_x}}\in \Delta], \quad k=0, 1, \cdots.$$
Then $p_k(x, A)$ is the $\P_x$ probability of the event that the process $Y$ goes to $A$
after exactly $k$ jumps from one set $C_y$ to another.
In the following, we assert that for $x\in \Delta,$
\begin{equation}\label{e:2.34}
\P_x(Y_{\tau_\Delta}\in A)=\sum_{k=0}^\infty p_k(x, A).
\end{equation}
 In fact, by the strong Markov property of $Y,$
$$\P_x(Y_{\tau_\Delta}\in A)=\P_x(Y_{\tau_{C_x}}\in A)+\E_x [\P_{Y_{\tau_{C_x}}}(Y_{\tau_\Delta}\in A); Y_{\tau_{C_x}}\in \Delta].$$
Define
$$r_0(x, A):=\P_x(Y_{\tau_\Delta}\in A), \quad r_{k+1}(x, A):=\E_x[r_k(Y_{\tau_{C_x}}, A); Y_{\tau_{C_x}}\in \Delta], \quad k=0, 1, \cdots.$$
Then $r_k(x, A)$ is the $\P_x$ probability of the event that the process $X$ goes to $A$
after more than $k$ jumps from one set $C_y$ to another.
By the induction argument,
\begin{equation}\label{e:2.35}
\P_x(Y_{\tau_\Delta}\in A)=\sum_{i=0}^k p_i(x, A)+r_{k+1}(x, A).
\end{equation}
 By Lemma \ref{L:2.6}, there exists $\rho_0=\rho_0(d, \lambda_0, \ell, \phi, R_0, \Lambda_0)\in (0, 1)$  such that  for any $j\geq 1,$
\begin{equation}\label{e:2.36}
 \P_x(Y_{\tau_{C_x}}\in \Delta)\leq \P_x(Y_{\tau_{B_x\cap D}}\in D)= 1-\P_x(Y_{\tau_{B_x\cap D}}\in D^c)\leq 1-\rho_0, \quad {\rm for} \quad x\in \Delta.
 \end{equation}
Thus by \eqref{e:2.36} and the induction argument,
$$r_{k+1}(x, A)\leq (1-\rho_0)^{k+1}\rightarrow 0  \quad {\rm as} \quad k\rightarrow\infty.$$
 This together with \eqref{e:2.35} establishes \eqref{e:2.34}.

 It follows from \cite[(1.5)]{CKS7} that there exists $c_1>0$ such that $j(u)\leq c_1|u|^{-(d+2)}$ for $u\in\R^d$ with $|u|\leq 1.$ Hence, by Proposition \ref{L:2.5} and the L\'evy system formula of $Y,$ there exist positive constant $c_k, k=2, 3$ such that
$$\begin{aligned}
\sup_{x\in \Delta}p_0(x, A)&= \sup_{x\in \Delta}\P_x(Y_{\tau_{C_x}}\in A)
=\sup_{x\in \Delta} \E_x\int_0^{\tau_{C_x}} \int_A J(Y_s, z)\,ds\,dz\\
&\leq c_2\sup_{x\in \Delta} \E_x\tau_{B_x} \int_\Upsilon^1 |z|^{-(d+2)}\,dz\\
&\leq c_32^{-2j}r^2\Upsilon^{-2}.
 \end{aligned}$$
  Let $k\geq 0.$  Suppose $\sup_{x\in \Delta}p_k(x, A)\leq c_3(1-\rho_0)^k2^{-2j}r^2 \Upsilon^{-2},$
  then by \eqref{e:2.36},
$$\begin{aligned}
\sup_{x\in \Delta}p_{k+1}(x, A)&=\sup_{x\in \Delta}\E_x[p_k(Y_{\tau_{C_x}}, A); Y_{\tau_{C_x}}\in \Delta]\\
&\leq c_3(1-\rho_0)^{k} 2^{-2j}r^2 \Upsilon^{-2}  \sup_{x\in \Delta}\P_x(Y_{\tau_{C_x}}\in \Delta)\\
&\leq c_3(1-\rho_0)^{k+1}2^{-2j}r^2 \Upsilon^{-2}.
\end{aligned}$$
Hence, by the induction and \eqref{e:2.34}, we have
$$
\sup_{x\in \Delta}\P_x(Y_{\tau_\Delta}\in A)\leq \sum_{k=0}^\infty\sup_{x\in \Delta}p_{k}(x, A)\leq \sum_{k=0}^\infty c_3(1-\rho_0)^{k}2^{-2j}r^2 \Upsilon^{-2} \leq \dfrac{c_4}{\rho_0}2^{-2j}r^2  \Upsilon^{-2}.
$$
\qed

\begin{lem}\label{L:2.8}
Suppose $d\geq 2$ and $H$ satisfies the assumption (A1).
Let $D$ be a $C^{1, \alpha}$ open set with characteristics $(R_0, \Lambda_0).$
There exists $C=C(d, \lambda_0, \ell, \phi, R_0, \Lambda_0)>0$ such that for any $z_0\in \partial D, r\in (0, (R_0\wedge a^{-1/2})/4)$ and $x\in D\cap B(z_0, r/4),$
$$\P_x(Y_{\tau_{D_r(z_0)}}\in  D_{2r}(z_0)\cap \Delta(z_0, r/2, 2r))\leq
C\P_x(Y_{\tau_{D_r(z_0)}}\in  D_{2r}(z_0)\setminus \Delta(z_0, r/2, 2r)),$$
where $D_r(z_0):=D\cap B(z_0, r).$
\end{lem}

\pf The proof mainly adapt the "box" method developed in \cite{Ba, BB1} in our case.
Without loss of generality, we assume the constant  $a=1$ in the assumption (A1).
Let $z_0\in \partial D$ and $r\in (0, R_0/4).$
Let
$$S:=D_{2r}(z_0)\cap \Delta(z_0, r/2, 2r), \quad U:= D_{2r}(z_0)\setminus \Delta(z_0, r/2, 2r).$$
Let $\omega_0(x):=\P_x(Y_{\tau_{D_r(z_0)}}\in  S)$ and $\omega_1(x):=\P_x(Y_{\tau_{D_r(z_0)}}\in  U).$

We define a decreasing sequence $\{r_i\}_{i\geq 0}$ by $r_0=r/2$ and
$$r_i:=\dfrac{r}{2}\left(1-\dfrac{3}{\pi^2}\sum_{j=1}^i \dfrac{1}{j^2}\right).$$
Note that $\sum_{j=1}^\infty j^{-2}=\pi^2/6.$ Hence, $r_i\in (r/4, r/2)$ for $i\geq 0.$
Define for each $j\geq 0,$
$$W_j:=\Delta(z_0, r2^{-(j+1)}, r_j)\setminus \Delta(z_0, r2^{-(j+2)}, r_j).$$
For each $j\geq 0,$ define
$$d_j:=\sup_{x\in \cup_{i=0}^j W_i}\dfrac{\omega_0(x)}{\omega_1(x)}.
$$
It is sufficient to show that there exists $C=C(d, \lambda_0, \ell, \phi, R_0, \Lambda_0)>0$ such that
$$\sup_{j \geq 0}d_j\leq C<\infty.$$

Note that for $x\in W_0,$ $\delta_D(x)>r/4.$
 By a similar argument in Lemma \ref{L:2.6'} with $U\cap \partial B(z_0, r)$ in place of $D\cap \partial B(z_0, r)$, we have there exits $c_1=c_1(d, \lambda_0, \ell, \phi, R_0, \Lambda_0)>0$ such that for any $r\in (0, R_0/4)$ and $x\in W_0,$
\begin{equation}\label{e:3.26}
w_1(x)=\P_x(Y_{\tau_{D_r(z_0)}}\in U)\geq \P_x(Y_{\tau_{D_r(z_0)}}\in U\cap \partial B(z_0, r))\geq c_1\dfrac{\delta_D(x)}{r}\geq c_1/4:=c_2.
\end{equation}
Hence, $w_0(x)\leq 1\leq c_2^{-1}w_1(x).$ Thus $d_0\leq c_2^{-1}.$

We define
$$J_j:=\cup_{i=0}^j W_i, \quad j\geq 0.$$
Let $\Omega:=\{Y_{\tau_{D_r(z_0)}}\in  S\}.$
Let $F_j:=\Delta(z_0, 2^{-(j+1)}r, (r_j+r_{j-1})/2).$
Let $\tau_j:=\tau_{F_j}.$
We have for $x\in W_j,$
\begin{equation}\label{e:3.24}\begin{aligned}
\P_x(\Omega)
&=\P_x(Y_{\tau_j}\in J_{j-1}; \Omega)+ \P_x(Y_{\tau_j}\in \Delta(z_0, r2^{-(j+1)}, r); \Omega)\\
&+\P_x(Y_{\tau_j}\in D_r(z_0)\setminus(J_{j-1}\cup \Delta(z_0, r2^{-(j+1)}, r)); \Omega)+\P_x(Y_{\tau_j}\in S).
\end{aligned}\end{equation}
By the strong Markov property of $Y$, for $x\in W_j,$
\begin{equation}\label{e:2.29}
\begin{aligned}
\P_x(Y_{\tau_j}\in J_{j-1}; \Omega)
&=\E_x (Y_{\tau_j}\in J_{j-1};  \P_{Y_{\tau_j}}(Y_{\tau_{D_r(z_0)}}\in S))\\
&\leq d_{j-1}\E_x (Y_{\tau_j}\in J_{j-1};  \P_{Y_{\tau_j}}(Y_{\tau_{D_r(z_0)}}\in U))\\
&\leq d_{j-1} \omega_1(x).
\end{aligned}
\end{equation}
Note that the distance between $W_j$ and $F_j^c$ is larger than $(r_{j-1}-r_j)/2.$
By applying Lemma  \ref{L:2.7} with $M=(r_{j-1}-r_j)/(r2^{-(j+2)})=\frac{3}{\pi^2}j^{-2}2^{(j+1)}$ and $s=r2^{-(j+1)}$,  for $x\in W_j,$
\begin{equation}\label{e:2.30}
\P_x(Y_{\tau_j}\in \Delta(z_0, r2^{-(j+1)}, r))\leq \P_x(Y_{\tau_{\Delta(x, s, Ms)}}\in \Delta(x, s, R_0/2))\leq c_3M^{-3}
\leq \pi^6 c_3j^{6}2^{-3j}/27.
\end{equation}
On the other hand, note that the distance between the points in $F_j$ and $D_r(z_0)\setminus(J_{j-1}\cup \Delta(z_0, r2^{-(j+1)}, r))$ is larger than
$(r_{j-1}-r_j)/2=3j^{-2}r/(2\pi^2),$ by Lemma \ref{L:2.15},  we have for $x\in W_j,$
\begin{equation}\label{e:2.31}
\P_x(Y_{\tau_j}\in D_r(z_0)\setminus(J_{j-1}\cup \Delta(z_0, r2^{-(j+1)}, r)); \Omega)\leq  c_4(2^{-j})^2 j^4.
\end{equation}
Similarly, since the distance between the points in $F_j$ and $S$ is larger than $r/2,$ by Lemma \ref{L:2.15},
\begin{equation}\label{e:2.32'}
\P_x(Y_{\tau_j}\in S)\leq  c_52^{-2j}, \quad x\in W_j.
\end{equation}

By Lemma \ref{L:2.6'}, there exists $c_6=c_6(d, \lambda_0, \ell, \phi, R_0, \Lambda_0)$ such that for any $r\in (0, R_0/4),$
\begin{equation}\label{e:3.18}\begin{aligned}
\omega_1(x)=\P_x(Y_{\tau_{D_r(z_0)}}\in U)
\geq c_6\delta_{D}(x)/r\geq c_62^{-j}, \quad x\in W_j.
\end{aligned}\end{equation}
Hence, by \eqref{e:3.24}-\eqref{e:3.18}, there exists $c_7$ such that for any  $r\in (0, R_0/4)$ and $x\in W_j,$
$$\begin{aligned}
w_0(x)=\P_x(\Omega)\leq d_{j-1} w_1(x)+c_72^{-2j}j^6\leq d_{j-1}w_1(x)+\dfrac{c_7}{c_6}w_1(x) 2^{-j}j^6.
\end{aligned}$$
Let $b_j:=\dfrac{c_7}{c_6}2^{-j}j^6.$ Then
$$d_j\leq d_{j-1}+b_j, \quad j\geq 1.$$
This implies that $\sup_{i\geq 1} d_i\leq d_0+\sum_{j=1}^\infty b_j<\infty.$
Thus we complete the proof.
\qed

\medskip
Let $D$ be a $C^{1, \alpha}$ domain with characteristics $(R_0, \Lambda_0)$ in $\R^d.$
Recall that for each $r\in (0, R_0/2),$ there exists $L=L(R_0, \Lambda_0, d)>1$
such that for any $z\in\partial D,$ there is a $C^{1,\alpha}$ connected open set  $U_{z, r}\subset D$  such that $D\cap B(z, r)\subset U_{z, r}\subset D\cap B(z, 2r)$ and  $r^{-1}U_{z, r}$  is a $C^{1, \alpha}$ open set with characteristics $(R_0/L, L\Lambda_0).$
Let $\delta_1$ be the constant in Proposition \ref{L:2.5} with $(R_0/L, L\Lambda_0)$ in place of $(R_0, \Lambda_0)$.

\begin{prp}\label{P:3.9}
Let $D$ be a $C^{1, \alpha}$ open set with characteristics $(R_0, \Lambda_0)$ in $\R^d$ with $d\geq 1$.
Suppose $H$ satisfies the assumption (A1).
There exists a positive constant $C=C(d, \lambda_0, \ell, \phi, R_0, \Lambda_0)$ such that for any $z_0\in \partial D, r\in (0, \delta_1/8)$ and
$x\in  D\cap B(z_0, r/4),$
\begin{equation}\label{e:3.30}
\P_x(Y_{\tau_{ D\cap B(z_0, r)}}\in D)\leq C\delta_{D}(x)/r.
\end{equation}
\end{prp}

\pf Without loss of generality, we assume the constant  $a=1$ in the assumption (A1).
 Let $z_0\in \partial D$ and $y_0$ be a point on $D\cap \partial B(z_0, 6r)$ with $\delta_{D}(y_0)>6\kappa r.$
For the simplicity of notation, let $D_r(z_0):=D\cap B(z_0, r).$
Let $A:=D_{2r}(z_0)\setminus \Delta(z_0, r/2, 2r)$ when $d\geq 2$ and let $A:=D_{2r}(z_0)\setminus D_{r/2}(z_0)$ when $d=1.$
Let $x\in D_{r/4}(z_0).$
By the strong Markov property of $Y$,
\begin{equation}\label{e:3.31}\begin{aligned}
\P_x(Y_{\tau_{D_{4r}(z_0)}}\in B(y_0, \kappa r))&\geq \E_x (\P_{Y_{\tau_{D_{r}(z_0)}}}(Y_{\tau_{D_{4r}(z_0)}}\in B(y_0, \kappa r)); Y_{\tau_{D_r(z_0)}}\in A)\\
&\geq \inf_{y\in A}\P_y(Y_{\tau_{D_{4r}(z_0)}}\in B(y_0, \kappa r)) \cdot \P_x(Y_{\tau_{D_r(z_0)}}\in A)\\
&\geq \inf_{y\in A}\P_y(Y_{\tau_{B(y, r/2)}}\in B(y_0, \kappa r)) \cdot \P_x(Y_{\tau_{D_r(z_0)}}\in A).
\end{aligned}\end{equation}
By the L\'evy system formula of $Y,$ Lemma \ref{L:2.5'} and Proposition \ref{L:2.5},
\begin{equation}\label{e:3.37}\begin{aligned}
\inf_{y\in A}\P_y(Y_{\tau_{B(y, r/2)}}\in B(y_0, \kappa r))
&=\inf_{y\in A}\E_y \int_0^{\tau_{B(y, r/2)}} \int_{B(y_0, \kappa r)} J(Y_s, z)\,ds\,dz\\
&\geq c_1\inf_{y\in A}\E_y \tau_{B(y, r/2)} \inf_{u\in B(y, r/2)}\int_{B(y_0, \kappa r)} j(|u-z|)\,dz\\
&\geq c_1r^{d+2} j(8r),
\end{aligned}\end{equation}
where the last inequality is due to that $j(\cdot)$ is decreasing and $|u-z|\leq 8r$ for $u\in B(y, r/2)$ and $z\in B(y_0, \kappa r).$

On the other hand, let $U_{z_0, 4r}$  be a $C^{1, \alpha}$ domain such that $D_{4r}(z_0)\subset U_{z_0, 4r}\subset D_{8r}(z_0)$ with characteristics $(4r R_0/L, L\Lambda_0/(4r)^\alpha).$ Then $r^{-1}U_{z_0, 4r}$  is a $C^{1, \alpha}$ open set with characteristics $(R_0/L, L\Lambda_0)$ and its diameter is less than $16.$
By the L\'evy system formula of $Y$, Proposition  \ref{L:2.5} and \eqref{e:2.45}, there exist  positive constants $c_k=c_k(d, \lambda_0, \ell, \phi, R_0, \Lambda_0), k=2, 3, 4$ such that for $r\in (0, \delta_1/8),$
\begin{equation}\label{e:3.32}\begin{aligned}
&\P_x(Y_{\tau_{D_{4r}(z_0)}}\in B(y_0, \kappa r))\leq c_2\E_x \tau_{D_{4r}(z_0)}\sup_{u\in D_{4r}(z_0)}\int_{B(y_0, \kappa r)} j(|u-z|)\,dz\\
&\leq c_2\E_x \tau_{U_{z_0, 4r}}\sup_{u\in D_{4r}(z_0)}\int_{B(y_0, \kappa r)} j(|u-z|)\,dz
\leq c_3\E_x \tau^W_{U_{z_0, 4r}}\cdot j(2r)r^d\leq c_4\delta_{D}(x)j(r)r^{d+1},
\end{aligned}\end{equation}
where the last inequality is due to that $j(\cdot)$ is decreasing and $|u-z|\geq r$ for $u\in D_{4r}(z_0)$ and $z\in B(y_0, \kappa r).$
Note that $j(r)$ is comparable to $j(8r)$ for $r\in (0, 1/8)$ by \eqref{e:2.9}. Hence, by \eqref{e:3.31}-\eqref{e:3.32}, there exists $c_5=c_5(d, \lambda_0, \ell, \phi, R_0, \Lambda_0)$ such that for $r\in (0, \delta_1/8)$ and $x\in D\cap B(z_0, r/4),$
\begin{equation}\label{e:3.40'}
\P_x(Y_{\tau_{D_r(z_0)}}\in A)\leq c_5\delta_{D}(x)/r.
\end{equation}
By Lemma \ref{L:2.8}, when $d\geq 2,$ there exists $c_6=c_6(d, \lambda_0, \ell, \phi, R_0, \Lambda_0)>1$ such that
$$\P_x(Y_{\tau_{D_r(z_0)}}\in D_{2r}(z_0))\leq c_6\P_x(Y_{\tau_{D_r(z_0)}}\in A).$$
When $d=1,$ it is easy to see that $$\P_x(Y_{\tau_{D_r(z_0)}}\in D_{2r}(z_0))=\P_x(Y_{\tau_{D_r(z_0)}}\in A).$$
Hence, it follows from \eqref{e:3.40'} that for $d\geq 1,$
\begin{equation}\label{e:3.39}
\P_x(Y_{\tau_{D_r(z_0)}}\in D_{2r}(z_0))\leq c_6\P_x(Y_{\tau_{D_r(z_0)}}\in A)\leq c_5c_6\delta_{D}(x)/r, \quad x\in D_{r/4}(z_0).
\end{equation}

By \eqref{e:2.9'} and \eqref{e:2.4}, there exists $c_7$ such that $j(z)\leq c_7|z|^{-(d+2)}$ for $|z|\leq 1$.
By the L\'evy system formula of $Y$, Proposition \ref{L:2.5} and \eqref{e:2.45},
\begin{equation}\label{e:3.40}\begin{aligned}
&P_x(Y_{\tau_{D_r(z_0)}}\in D\setminus D_{2r}(z_0))\\
&=\E_x \int_0^{\tau_{D_r(z_0)}}\int_{D\setminus D_{2r}(z_0)} J(Y_s, z)\,ds\,dz\\
&\leq c_8\E_x \tau_{D_r(z_0)}\sup_{u\in D_r(z_0)}\int_{D\setminus D_{2r}(z_0)}j(|u-z|)\,dz\\
&\leq c_9\E_x \tau^W_{D_r(z_0)}\left(\int_{r\leq |z|\leq 1}|z|^{-(d+2)}\,dz+\int_{|z|>1}j(|z|)\,dz\right)\\
&\leq c_{10}\dfrac{\delta_D(x)}{r}.
\end{aligned}\end{equation}
By combining \eqref{e:3.39} and \eqref{e:3.40}, the desired conclusion is obtained.
\qed

\section{Upper and Lower bound estimates}
\subsection{Upper bound estimates}
In this section, we shall establish the upper bound estimates of the Dirichlet heat kernel of $Y$ in a $C^{1, \alpha}$ open set $D$ in Theorem \ref{T1}.
By a very similar argument in \cite[Lemma 3.1]{CKS7} and the strong Markov property of $Y,$ we have the following Lemma.

\begin{lem}\label{L:4.1}
Suppose that $U_1, U_3, E$ are open subsets of $\R^d$ with $U_1, U_3\subset E$ and ${\rm dist}(U_1, U_3)>0.$
Let $U_2:=E\setminus (U_1\cup U_3).$ If $x\in U_1$ and $y\in U_3,$ then for every $t>0,$
\begin{eqnarray}
p_E(t, x, y)
&\leq &\E_x(p_E(t-\tau_{U_1}, Y_{\tau_{U_1}}, y); Y_{\tau_{U_1}}\in U_2, \tau_{U_1}<t)\nonumber\\
&\quad &+\E_x(p_E(t-\tau_{U_1}, Y_{\tau_{U_1}}, y); Y_{\tau_{U_1}}\in U_3; \tau_{U_1}<t)\nonumber\\
&\leq & \P_x(Y_{\tau_{U_1}}\in U_2)\left(\sup_{s<t, z\in U_2}p_E(s, z, y)\right)\nonumber\\
&+& \int_0^t \P_x(\tau_{U_1}>s)\P_y(\tau_E>t-s)\,ds \left(\sup_{u\in U_1, z\in U_3} J(u, z)\right)\label{e:4.0}\\
&\leq& \P_x(Y_{\tau_{U_1}}\in U_2)\left(\sup_{s<t, z\in U_2}p_E(s, z, y)\right)+(t\wedge \E_x\tau_{U_1})\left(\sup_{u\in U_1, z\in U_3} J(u, z)\right)\label{e:4.0'}.
\end{eqnarray}
\end{lem}

Recall that $\delta_1\in (0, R_0)$ is the constant in Proposition \ref{L:2.5}.

\begin{prp}\label{P:4.2}
Suppose that $D$ is a $C^{1, \alpha}$ open set in $\R^d$ with characteristics $(R_0, \Lambda_0).$
If $D$ is bounded,  assume that $H$ satisfies the assumption (A1).
If $D$ is unbounded,  assume that $H$ satisfies the assumptions (A1) and (A2).
For every $T>0,$ there exist positive constants $C=C(d, \lambda_0, \ell, \phi, R_0, \Lambda_0, T), a_U=a_U(d, \lambda_0, \ell, \phi)$ and $b_U=b_U(d, \lambda_0, \ell, \phi)$ such that for  all $x, y\in D$ and $t\in (0, T),$
$$\begin{aligned}
&p_{D}(t, x, y)\leq C\left(1\wedge \dfrac{\delta_{D}(x)}{\sqrt t}\right)\\
&\quad\times\left[t^{-d/2}\wedge\left(t^{-d/2}e^{-|x-y|^2/(4b_Ut)}+\dfrac{tH(|x-y|^{-2})}{|x-y|^d}
+\phi^{-1}(1/t)^{d/2}e^{-\frac{a_U}{4}|x-y|^2\phi^{-1}(1/t)}\right)\right].
\end{aligned} $$
\end{prp}

\pf Let $x\in D.$ In view of \eqref{e:2.26'} and \eqref{e:2.29'}, we only need to prove the theorem for $\delta_{D}(x)
\leq \sqrt{t}\delta_1/(16\sqrt T).$
Let  $z_0\in \partial D$ be such that $\delta_{D}(x)=|x-z_0|.$ Let $t\in (0, T).$ Let
$$U_1:= B(z_0, \sqrt{t} \delta_1/(8\sqrt T))\cap D.$$
Let $\tilde U_1$  be a $C^{1, \alpha}$ domain such that $U_1\subset \tilde U_1\subset B(z_0, \sqrt{t}\delta_1/(4\sqrt T))\cap D$ and
$\frac{8\sqrt T}{\delta_1 \sqrt{t}}\tilde U_1$  is a $C^{1, \alpha}$ open set with characteristics $(R_0/L, L\Lambda_0).$
By Proposition \ref{L:2.5}, there exists $c_1=c_1(d, \lambda_0, \ell, \phi, R_0, \Lambda_0)$ such that
\begin{equation}\label{e:4.1'}
\E_x \tau_{U_1}\leq \E_x \tau_{\tilde U_1}\leq c_1\delta_D(x)\sqrt{t}.
\end{equation}
By Proposition \ref{P:3.9} with $r=\sqrt{t} \delta_1/(8\sqrt T),$ there exists $c_1=c_1(d, \lambda_0, \ell, \phi, R_0, \Lambda_0)>0$ such that
\begin{equation}\label{e:4.4n}
\P_x(Y_{\tau_{U_1}}\in D)\leq c_1\left(1\wedge \frac{\delta_{D}(x)}{\sqrt t}\right).
\end{equation}
By applying the strong Markov property of $Y,$ \eqref{e:4.1'} and \eqref{e:4.4n}, there exists $c_2=c_2(d, \lambda_0, \ell, \phi, R_0, \Lambda_0)>0$ such that for any $t\in (0, T),$
\begin{equation}\label{e:4.1}\begin{aligned}
\P_x(\tau_{D}>t)
&\leq \P_x(\tau_{U_1}>t)+\P_x(Y_{\tau_{U_1}}\in D\setminus U_1)\\
&\leq 1\wedge \dfrac{\E_x \tau_{U_1}}{t}+\P_x(Y_{\tau_{U_1}}\in D\setminus U_1)\\
&\leq c_2\left(1\wedge \frac{\delta_{D}(x)}{\sqrt t}\right).
\end{aligned}\end{equation}

Note that by Theorem \ref{T0}, there exist positive constants $c_3>0, b_U=b_U(d, \lambda_0, \phi)$ and $a_U=a_U(d, \lambda_0, \phi)$ such that
\begin{equation}\label{e:4.6n}
p(t, x, y)\leq c_3t^{-d/2}\wedge\left(t^{-d/2}e^{-|x-y|^2/(b_Ut)}+\dfrac{tH(|x-y|^{-2})}{|x-y|^d}
+\phi^{-1}(1/t)^{d/2}e^{-a_U|x-y|^2\phi^{-1}(1/t)}\right)
\end{equation}
 holds for $x, y\in D$ and $t\in (0, T)$ when $D$ is  bounded  under the assumption (A1), and holds for $x, y\in\R^d$ and $t>0$  under the assumption (A2).

Now we deal with two cases separately.
Let $c_0:=(d/2)\vee [(dc_L^{-1/\gamma}T^{1/\gamma-1}\phi^{-1}(1)^{-1})/(2b_Ua_U)]\vee [\delta_1^2/(4b_UT)]$ and let $x, y\in D.$

Case 1: $|x-y|\leq 2(b_Uc_0)^{1/2}\sqrt t.$ By the semigroup property, \eqref{e:4.1} and \eqref{e:4.6n},
$$\begin{aligned}
p_{D}(t, x, y)&=\int_{D} p_{D}(t/2, x, z)p_{D}(t/2, z, y)\,dz\\
&\leq \sup_{z, w\in D} p_{D}(t/2, z, w)\int_{D} p_{D}(t/2, x, z)\,dz\\
&\leq c_4 t^{-d/2}\P_x(\tau_{D}>t/2)\\
&\leq c_5t^{-d/2}\left(1\wedge \frac{\delta_{D}(x)}{\sqrt t}\right).
\end{aligned}$$
Since $|x-y|^2/(4b_U\sqrt{t})\leq c_0,$ we have
\begin{equation}\label{e:4.6}
p_{D}(t, x, y)\leq c_6t^{-d/2}e^{-|x-y|^2/(4b_Ut)}\left(1\wedge \frac{\delta_{D}(x)}{\sqrt t}\right).
\end{equation}

Case 2: $|x-y|>2(b_Uc_0)^{1/2}\sqrt t.$ Let
$$U_3:=\{z\in D: |z-x|\geq |x-y|/2\}, \quad U_2:=D\setminus (U_1\cup U_3).$$
Note that $|x-y|>2(b_Uc_0)^{1/2}\sqrt t\geq \delta_1\sqrt {t}/\sqrt{T}.$
Hence the distance between $U_1$ and $U_3$ is larger than $\delta_1\sqrt {t}/(4\sqrt{T}).$
For $z\in U_2,$ it is easy to see that
$$\dfrac{|x-y|}{2}\leq |z-y|\leq \dfrac{3|x-y|}{2}.$$
Recall that $p^0(t, x, y)$ is the transition density function of $Y^0=W_{S_t}.$
By\eqref{e:4.6n}, \cite[(3.23)]{BK} and the choice of $c_0$, we have
\begin{equation}\label{e:4.7}
\begin{aligned}
&\sup_{s\leq t, z\in U_2} p(s, z, y)\leq c_7\sup_{s\leq t, z\in U_2} p^0(s, z, y)\\
&\leq c_8\left(t^{-d/2}e^{-|x-y|^2/(4b_Ut)}+\dfrac{tH(|x-y|^{-2})}{|x-y|^d}+\phi^{-1}(1/t)^{d/2}e^{-\frac{a_U}{4}|x-y|^2\phi^{-1}(1/t)}\right).
\end{aligned}
\end{equation}

For $u\in U_1$ and $z\in U_3,$ we have $|z-x|>|x-y|/2,$ thus by the choice of $c_0,$
$$|u-z|\geq |z-x|-|x-u|\geq |z-x|/2\geq |x-y|/4.$$
By Lemma 2.3 in \cite{BK} (cf. \cite[Lemma 2.1]{M}),
\begin{equation}\label{e:4.8}
H(\lambda s)\leq \lambda^2 H(s), \quad \lambda\geq 1, s>0.
\end{equation}
Hence, by \eqref{e:2.17},
\begin{equation}\label{e:4.5'}
\sup_{u\in U_1, z\in U_3}J(u, z)\leq c_9\sup_{u\in U_1, z\in U_3}\dfrac{H(|u-z|^{-2})}{|u-z|^d}\leq \dfrac{c_9 4^{d+4}H(|x-y|^{-2})}{|x-y|^d}.
\end{equation}
Consequently, by Lemma \ref{L:4.1}, \eqref{e:4.1'}-\eqref{e:4.4n}, \eqref{e:4.7} and \eqref{e:4.5'}, for $|x-y|>2(b_Uc_0)^{1/2}\sqrt t,$
$$\begin{aligned}
&p_{D}(t, x, y)\\
\leq &\P_x(Y_{\tau_{U_1}}\in U_2)\left(\sup_{s<t, z\in U_2}p(s, z, y)\right)+(t\wedge \E_x\tau_{U_1})\left(\sup_{u\in U_1, z\in U_3} J(u, z)\right)\\
\leq & c_{10}\left(1\wedge \dfrac{\delta_{D}(x)}{\sqrt t}\right)\left[t^{-d/2}\wedge\left(t^{-d/2}e^{-|x-y|^2/(4b_Ut)}+\dfrac{tH(|x-y|^{-2})}{|x-y|^d}
+\phi^{-1}(1/t)^{d/2}e^{-\frac{a_U}{4}|x-y|^2\phi^{-1}(1/t)}\right)\right].
\end{aligned}$$
The proof is complete.
\qed

\medskip

{\bf Proof of Theorem \ref{T1} (i).} Fix $T>0.$ Let $t\in (0, T]$ and $x, y\in D.$
By Proposition \ref{P:4.2}, Theorem \ref{T0} and the symmetry of $p_D(t, x, y)$ in $(x, y)$, we only need to prove Theorem \ref{T1} (i) when $\delta_D(x)\vee\delta_D(y)\leq \delta_1\sqrt{t}/(16\sqrt{T})\leq \delta_1/16.$
The proof is along the line of the proof of Proposition \ref{P:4.2}.
Define $U_1$ in the same way as in the proof of Proposition \ref{P:4.2}.
Let $a_U$ and $b_U$ be the constants in \eqref{e:4.6n}.
Let $c_1:=((d+1)/2)\vee [(dc_L^{-1/\gamma}T^{1/\gamma-1}\phi^{-1}(1)^{-1})/(a_Ub_U)].$
We estimate $p_D(t, x, y)$ by considering the following two cases.

Case 1: $|x-y|\leq 4(b_Uc_1)^{1/2}\sqrt t.$ By the semigroup property, Proposition \ref{P:4.2} and \eqref{e:4.1},
$$\begin{aligned}
p_{D}(t, x, y)&=\int_{D} p_{D}(t/2, x, z)p_{D}(t/2, z, y)\,dz\\
&\leq \sup_{z\in D} p_{D}(t/2, y, z)\int_{D} p_{D}(t/2, x, z)\,dz\\
&\leq c_2 t^{-d/2}\left(1\wedge \frac{\delta_{D}(y)}{\sqrt t}\right)\P_x(\tau_{D}>t/2)\\
&\leq c_3t^{-d/2}\left(1\wedge \frac{\delta_{D}(x)}{\sqrt t}\right)\left(1\wedge \frac{\delta_{D}(y)}{\sqrt t}\right)\\
&\leq c_4t^{-d/2}e^{-|x-y|^2/(16b_Ut)}\left(1\wedge \frac{\delta_{D}(x)}{\sqrt t}\right)\left(1\wedge \frac{\delta_{D}(y)}{\sqrt t}\right).
\end{aligned}$$

Case 2: $|x-y|>4(b_Uc_1)^{1/2}\sqrt t.$ Define $U_2$ and $U_3$ in the same way as in the proof of Proposition \ref{P:4.2}.
Note that for $z\in U_2,$ $\dfrac{|x-y|}{2}\leq |z-y|\leq \dfrac{3|x-y|}{2}.$
By Proposition \ref{P:4.2}, the choice of $c_1$ and a very similar argument in \cite[(3.29)-(3.30)]{BK},
\begin{equation}\label{e:4.4'}
\begin{aligned}
&\sup_{s\leq t, z\in U_2}p_{D}(s, z, y)\\
&\leq c_5\left(1\wedge \dfrac{\delta_{D}(y)}{\sqrt{t}}\right)\left[t^{-d/2}\wedge\left(t^{-d/2}e^{-|x-y|^2/(4b_Ut)}+\dfrac{tH(|x-y|^{-2})}{|x-y|^d}
+\phi^{-1}(1/t)^{d/2}e^{-\frac{a_U}{4}|x-y|^2\phi^{-1}(1/t)}\right)\right].
\end{aligned}\end{equation}
On the other hand, by \eqref{e:4.1},
\begin{equation}\label{e:4.5}
\begin{aligned}
&\int_0^t \P_x(\tau_{U_1}>s)\P_y(\tau_D>t-s)\,ds\leq \int_0^t \P_x(\tau_D>s)\P_y(\tau_D>t-s)\,ds\\
&\leq c_6 \int_0^t \dfrac{\delta_D(x)}{\sqrt{s}}\dfrac{\delta_D(y)}{\sqrt{t-s}}\,ds=c_6\delta_D(x)\delta_D(y)\int_0^1 \dfrac{1}{\sqrt{r(1-r)}}
=c_7\delta_D(x)\delta_D(y).
\end{aligned}\end{equation}

Thus by \eqref{e:4.0} together with \eqref{e:4.5'} and \eqref{e:4.4'}-\eqref{e:4.5},
$$\begin{aligned}
&p_{D}(t, x, y)\\
\leq & \P_x(Y_{\tau_{U_1}}\in U_2)\left(\sup_{s<t, z\in U_2}p_D(s, z, y)\right)\\
&+ \int_0^t \P_x(\tau_{U_1}>s)\P_y(\tau_D>t-s)\,ds \left(\sup_{u\in U_1, z\in U_3} J(u, z)\right)\\
\leq & c_8\left(1\wedge \dfrac{\delta_{D}(x)}{\sqrt t}\right)\left(1\wedge \dfrac{\delta_{D}(y)}{\sqrt t}\right) \\ &\quad \times\left[t^{-d/2}\wedge\left(t^{-d/2}e^{-|x-y|^2/(4b_Ut)}+\dfrac{tH(|x-y|^{-2})}{|x-y|^d}
+\phi^{-1}(1/t)^{d/2}e^{-\frac{a_U}{4}|x-y|^2\phi^{-1}(1/t)}\right)\right].
\end{aligned}$$
The proof is complete.
\qed

\subsection{Lower bound estimates}

In this section, let $D$ be a $C^{1, \alpha}$ open set in $\R^d$ with characteristics $(R_0, \Lambda_0).$  We shall establish the lower bound estimates of the Dirichlet heat kernel of $Y$ in $D$  in Theorem \ref{T1}.

By the result in \cite{Cho}, for fixed $T>0,$ there exist positive constants $c_k=c_k(d, \lambda_0, \ell, R_0, \Lambda_0, T), k=1, 2$ such that for any $ x, y\in D$ and $t\in (0, T),$
\begin{equation}\label{e:4.14}
p^X_D(t, x, y)\geq c_1 \left(1\wedge \dfrac{\delta_D(x) \delta_D(y)}{t}\right)t^{-d/2} \exp(-c_2\dfrac{|x-y|^2}{t}).
\end{equation}
By \eqref{e:4.14}, \cite[Lemma 2.4]{BK} and a similar argument in \cite[Proposition 3.1]{BK} with $X$ in place of Brownian motion $W$, we have the following Lemma.

\begin{lem}\label{L:4.3}
Suppose that $D$ is a $C^{1, \alpha}$ open set in $\R^d$ with characteristics $(R_0, \Lambda_0).$ If $D$ is bounded, we assume that $H$ satisfies $L^a(\gamma, C_L)$ for some $a>0.$ If $D$ is unbounded, we assume that $H$ satisfies $L^0(\gamma, C_L)$ and the path distance in each connected component of $D$ is comparable to the Euclidean distance with characteristic $\chi_1.$ For each $T>0,$ there exist positive constants $C_1=C_1(d, \lambda_0, \ell, \phi, \chi_1, R_0, \Lambda_0, T)$ and $C_2=C_2(d, \lambda_0, \phi, \ell,  \chi_1, R_0, \Lambda_0)$ such that for all $ t\in (0, T]$ and $x, y\in D$  in the same connected component of $D,$
$$p_{D}(t, x, y)\geq C_1\left(1\wedge \dfrac{\delta_{D}(x)}{\sqrt t}\right)\left(1\wedge \dfrac{\delta_{D}(y)}{\sqrt t}\right)\phi^{-1}(t^{-1})^{d/2}e^{-C_2|x-y|^2\phi^{-1}(t^{-1})}.$$
\end{lem}

Denote by $T_t$ the subordinator with the Laplace exponent $\phi.$ Then $S_t=t+T_t$. By \cite[Proposition III.8]{Ber} and the Markov property of $T_t,$ for each $b>0$ and $T>0,$
there exists $c=c(\phi, b, T)$ such that
$$\P(T_t\leq bt)\geq c, \quad t\leq T.$$
Hence, by \eqref{e:4.14} and a similar argument in \cite[Lemma 2.1]{CKS7} with $X$ in place of $W$, we have the following Lemma.

\begin{lem}\label{L:4.4}
Suppose that $D$ is a $C^{1, \alpha}$ open set in $\R^d$ with characteristics $(R_0, \Lambda_0)$  and the path distance in each connected component of $D$ is comparable to the Euclidean distance with characteristic $\chi_1.$ For each $T>0,$ there exist positive constants $C_1=C_1(d, \lambda_0, \ell, \phi, \chi_1, R_0, \Lambda_0, T)$ and $C_2=C_2(d, \lambda_0, \ell,  \chi_1, R_0, \Lambda_0)$ such that for all $ t\in (0, T]$ and $x, y\in D$  in the same connected component of $D,$
$$p_{D}(t, x, y)\geq C_1\left(1\wedge \dfrac{\delta_{D}(x)}{\sqrt t}\right)\left(1\wedge \dfrac{\delta_{D}(y)}{\sqrt t}\right)t^{-d/2}e^{-C_2|x-y|^2/t}.$$
\end{lem}

The following two Lemmas can be obtained by Lemmas \ref{L:4.3}-\ref{L:4.4} and the same argument as \cite[Lemmas 3.2-3.3]{BK}.
We omit the proof here.

\begin{lem}\label{L:4.5}
For each positive constant $\varrho,$ there exists $c=c(d, \lambda_0, \ell, \varrho, \phi)>0$ such that for all $x\in\R^d$ and $r>0,$
$$\inf_{y\in B(x, r)}\P_y(\tau_{B(x, 2r)}\geq \varrho r^2)\geq c.$$
\end{lem}

\begin{lem}\label{L:4.6}
Suppose $H$ satisfies $L^a(\gamma, c_L)$ and $U^a(\delta, C_U)$ with  $\delta\leq 1$ for some $a>0$ ($L^0(\gamma, c_L)$ and $U^0(\delta, C_U)$, respectively).
Then for each $T>0, M>0$ and $b>0,$
 there exists $c=c(b, \phi)>0$ such that for all $t\in (0, T)$ and $u, v\in\R^d$ with $|u-v|\leq M/2$ ($u, v\in\R^d$, respectively)
$$p_{E}(t, u, v)\geq c(t^{-d/2}\wedge t|u-v|^{-d}H(|u-v|^{-2})),$$
where $E:=B(u, bt^{1/2})\cup B(v, bt^{1/2}).$
\end{lem}

{\bf Proof of Theorem \ref{T1} (ii)}
 Let  $D$ be a $C^{1, \alpha}$ open set with characteristics $(R_0, \Lambda_0)$ in $\R^d$. Then $D$ is a Lipschitz open set with characteristics $(R_0, \Lambda_0)$ in $\R^d$.  By the "corkscrew condition" for Lipschitz domain (see e.g. \cite[Lemma 6.6]{ChZ}), there exist constants $C_0\geq 1$ and $r_0\in (0, R_0)$ such that for any $z\in \partial D$ and $0<r\leq r_0,$ we can find a point $A=A_r(z)$ in $D$ satisfying $|A-z|\leq C_0r$ and $\delta_{D}(A)\geq r$.
Set $T_0=(r_0/4)^2.$ By considering the cases $\delta_D(x)<r_0$ and $\delta_D(x)>r_0,$  there exists $L_0=L(r_0)>1$ such that for any $t\in (0, T_0]$ and $x, y\in D,$ one can choose $A^t_x\in D\cap B(x, L_0 \sqrt{t})$ and $A^t_y\in D\cap B(y, L_0 \sqrt{t})$ so that $B(A^t_x, 2\sqrt{t})$ and $B(A^t_y, 2\sqrt{t})$are subsets  of the connected components of $D$ that contains $x$ and $y$ respectively.

We first consider the case $t\in (0, T_0].$
Note that for $u\in B(A^t_x, \sqrt{t}),$
$$\delta_{D}(u)\geq \sqrt{t} \quad \mbox{and} \quad |x-u|\leq |x-A^t_x|+|A^t_x-u|\leq (L_0+1)\sqrt{t}.$$
Then by Lemma \ref{L:4.4}, for $t\in (0, T_0],$
\begin{equation}\label{e:4.9}
\begin{aligned}
\int_{B(A^t_x, \sqrt{t})}p_{D}(t/3, x, u)\,du
&\geq c_1\left(1\wedge \dfrac{\delta_{D}(x)}{\sqrt t}\right)\int_{B(A^t_x, \sqrt{t})}\left(1\wedge \dfrac{\delta_{D}(u)}{\sqrt t}\right)t^{-d/2}e^{-c_2|x-u|^2/t}\,du\\
&\geq c_1\left(1\wedge \dfrac{\delta_{D}(x)}{\sqrt t}\right)t^{-d/2}e^{-c_2(L_0+1)^2}|B(A^t_x, \sqrt{t})|
\geq c_3\left(1\wedge \dfrac{\delta_{D}(x)}{\sqrt{ t}}\right).
\end{aligned}\end{equation}
Similarly, for $t\in (0, T_0],$
\begin{equation}\label{e:4.10}
\int_{B(A^t_y, \sqrt{t})}p_{D}(t/3, y, u)\,du\geq c_4\left(1\wedge \dfrac{\delta_{D}(y)}{\sqrt{ t}}\right).
\end{equation}
By the semigroup property, for $t\in (0, T_0],$
\begin{equation}\label{e:4.4}
p_{D}(t, x, y)\geq \int_{B(A^t_x, \sqrt{t})}\int_{B(A^t_y, \sqrt{t})}p_{D}(t/3, x, u)p_{D}(t/3, u, v)p_{D}(t/3, v, y)\,du\,dv.
\end{equation}
In the following,
we consider the cases $|x-y|\geq \sqrt{t}/8$ and $|x-y|< \sqrt{t}/8$ separately.

Case 1: Suppose $|x-y|\geq \sqrt{t}/8$ and $t\in (0, T_0).$ By \eqref{e:4.9}-\eqref{e:4.10} and Lemma \ref{L:4.6},
\begin{equation}\label{e:4.12}\begin{aligned}
&p_{D}(t, x, y)\\
\geq &\int_{B(A^t_x, \sqrt{t})}\int_{B(A^t_y, \sqrt{t})}p_D(t/3, x, u)p_{B(A^t_x, \sqrt{t})\cup B(A^t_y, \sqrt{t})}(t/3, u, v)p_D(t/3, v, y)\,du\,dv\\
\geq &c_5\left(1\wedge \dfrac{\delta_{D}(x)}{\sqrt{ t}}\right)\left(1\wedge \dfrac{\delta_{D}(y)}{\sqrt{ t}}\right)
\inf_{(u, v)\in B(A^t_x, \sqrt{t})\times B(A^t_y, \sqrt{t})}\left(t^{-d/2}\wedge \left(t\dfrac{H(|u-v|^{-2})}{|u-v|^d}\right)\right)\\
\geq &c_6\left(1\wedge \dfrac{\delta_{D}(x)}{\sqrt{ t}}\right)\left(1\wedge \dfrac{\delta_{D}(y)}{\sqrt{ t}}\right)
\left(t^{-d/2}\wedge \left(t\dfrac{H(|x-y|^{-2})}{|x-y|^d}\right)\right),
\end{aligned}\end{equation}
where in the last inequality, we used $|u-v|\leq c|x-y|$ for $|x-y|\geq \sqrt{t}/8$ and \eqref{e:4.8}.
By Lemma 3.4 in \cite{BK}, for any given positive constants $c_7, c_8, R$ and $T,$ there is a positive constant $c_9$ such that
\begin{equation}\label{e:4.19}
t^{-d/2}e^{-r^2/(4c_7t)}+\phi^{-1}(1/t)^{d/2}e^{-c_8r^2\phi^{-1}(1/t)}\leq c_9tH(r^{-2})r^{-d}
\end{equation}
for any $r\geq R$ and $t\in (0, T).$
By combining \eqref{e:4.12}-\eqref{e:4.19}, Lemmas \ref{L:4.3}-\ref{L:4.4} and by considering two cases when  $x$ and $y$ are contained in a connected component of $D$ or in two distinct components of $D$ separately, we obtain
$$\begin{aligned}
&p_{D}(t, x, y)\\
\geq & c_{10}\left(1\wedge \dfrac{\delta_{D}(x)}{\sqrt t}\right)\left(1\wedge \dfrac{\delta_{D}(y)}{\sqrt{ t}}\right)\\
&\quad \times\left[t^{-d/2}\wedge\left(t^{-d/2}e^{-|x-y|^2/(4c_{11}t)}+\dfrac{tH(|x-y|^{-2})}{|x-y|^d}
+\phi^{-1}(1/t)^{d/2}e^{-c_{12}|x-y|^2\phi^{-1}(1/t)}\right)\right].
\end{aligned}$$

Case 2: Suppose $|x-y|<\sqrt{t}/8$ and $t\in (0, T_0).$ In this case $x$ and $y$ are in the same connected component of $D.$
Then for $(u, v)\in B(A^t_x, \sqrt{t})\times B(A^t_y, \sqrt{t}),$
$$|u-v|\leq 2(1+L_0) \sqrt{t}+|x-y|\leq 2((1+L_0)+1/8)\sqrt{t}.$$
Then by Lemma \ref{L:4.4}, for $(u, v)\in B(A^t_x, \sqrt{t})\times B(A^t_y, \sqrt{t}),$
$$p_{D}(t/3, u, v)\geq c_{13}\left(1\wedge \dfrac{\delta_{D}(u)}{\sqrt t}\right)\left(1\wedge \dfrac{\delta_{D}(v)}{\sqrt t}\right)t^{-d/2}e^{-c_{14}|u-v|^2/t}\geq c_{15}t^{-d/2}.$$
Hence, by applying \eqref{e:4.9}-\eqref{e:4.4}, for $t\in (0, T_0],$
$$\begin{aligned}
&p_{D}(t,x,y)\\
&\geq c_{16}\left(1\wedge \dfrac{\delta_{D}(x)}{\sqrt{ t}}\right)\left(1\wedge \dfrac{\delta_{D}(y)}{\sqrt{ t}}\right)t^{-d/2}\\
&\geq  c_{16}\left(1\wedge \dfrac{\delta_{D}(x)}{\sqrt t}\right)\left(1\wedge \dfrac{\delta_{D}(y)}{\sqrt{ t}}\right)\\
&\quad \times\left[t^{-d/2}\wedge\left(t^{-d/2}e^{-|x-y|^2/(c_{13}t)}+\dfrac{tH(|x-y|^{-2})}{|x-y|^d}
+\phi^{-1}(1/t)^{d/2}e^{-c_{14}|x-y|^2\phi^{-1}(1/t)}\right)\right].
\end{aligned}$$

When $T>T_0$ and $t\in (T_0, T],$ observe that $T_0/3\leq t-2T_0/3\leq (T/T_0-2/3)T_0,$ that is, $t-2T_0/3$ is comparable to $T_0/3$ with some universal constants that depend only on $T$ and $T_0.$ Using the inequality
$$p_{D}(t, x, y)\geq \int_{B(A^{T_0}_x, \sqrt{T_0})}\int_{B(A^{T_0}_y, \sqrt{T_0})}p_{D}(T_0/3, x, u)p_D(t-2T_0/3, u, v)p_{D}(T_0/3, v, y)\,du\,dv$$
instead of \eqref{e:4.4} and by considering the case $|x-y|\geq \sqrt{T_0}/8$ and $|x-y|< \sqrt{T_0}/8$ separately, we obtain by the same argument as above that the lower bound holds for $t\in (T_0, T]$ and hence for $t\in (0, T].$
\qed

\medskip

{\bf Proof of Theorem \ref{T1}(iii) and Corollary \ref{C1}} The proof of  Theorem \ref{T1}(iii) is  the same as \cite[Theorem \ref{T1}(iii)-(iv)]{CKS7}. Corollary \ref{C1} follows by Theorem \ref{T1} (i)-(ii) and the same argument in \cite[Corollary 1.4]{BK}.  We omit the proof here.
\qed

\vskip 0.3truein

{\bf Jie-Ming Wang}

School of Mathematics and Statistics,  Beijing Institute of Technology,
Beijing 100081, P. R. China

E-mail: \texttt{wangjm@bit.edu.cn}

\end{document}